

\catcode`\@=11 
 
\def\nolabels{\def\wrlabel##1{}\def\eqlabel##1{}\def\reflabel##1{}}
\def\writelabels{\def\wrlabel##1{\leavevmode\vadjust{\rlap{\smash%
{\line{{\escapechar=` \hfill\rlap{\sevenrm\hskip.03in\string##1}}}}}}}%
\def\eqlabel##1{{\escapechar-1\rlap{\sevenrm\hskip.05in\string##1}}}%
\def\thlabel##1{{\escapechar-1\rlap{\sevenrm\hskip.05in\string##1}}}%
\def\reflabel##1{\noexpand\llap{\noexpand\sevenrm\string\string\string##1}}}
\nolabels
\global\newcount\secno \global\secno=0
\global\newcount\meqno \global\meqno=1
\global\newcount\mthno \global\mthno=1
\global\newcount\mexno \global\mexno=1
\global\newcount\mquno \global\mquno=1
\global\newcount\tblno \global\tblno=1
\def\newsec#1{\global\advance\secno by1 
\global\subsecno=0\xdef\secsym{\the\secno.}\global\meqno=1\global\mthno=1
\global\mexno=1\global\mquno=1\global\figno=1\global\tblno=1

\bigbreak\medskip\noindent{\bf\the\secno. #1}\writetoca{{\secsym} {#1}}
\par\nobreak\medskip\nobreak}
\xdef\secsym{}
\global\newcount\subsecno \global\subsecno=0
\def\subsec#1{\global\advance\subsecno by1 \global\subsubsecno=0
\xdef\subsecsym{\the\subsecno.}
\bigbreak\noindent{\bf\secsym\the\subsecno. #1}\writetoca{\string\quad
{\secsym\the\subsecno.} {#1}}\par\nobreak\medskip\nobreak}
\xdef\subsecsym{}
\global\newcount\subsubsecno \global\subsubsecno=0
\def\subsubsec#1{\global\advance\subsubsecno by1
\bigbreak\noindent{\it\secsym\the\subsecno.\the\subsubsecno.
                                   #1}\writetoca{\string\quad
{\the\secno.\the\subsecno.\the\subsubsecno.} {#1}}\par\nobreak\medskip\nobreak}
\global\newcount\appsubsecno \global\appsubsecno=0
\def\appsubsec#1{\global\advance\appsubsecno by1 \global\subsubsecno=0
\xdef\appsubsecsym{\the\appsubsecno.}
\bigbreak\noindent{\it\secsym\the\appsubsecno. #1}\writetoca{\string\quad
{\secsym\the\appsubsecno.} {#1}}\par\nobreak\medskip\nobreak}
\xdef\appsubsecsym{}
\def\appendix#1#2{\global\meqno=1\global\mthno=1\global\mexno=1
\global\figno=1\global\tblno=1
\global\subsecno=0\global\subsubsecno=0
\global\appsubsecno=0
\xdef\appname{#1}
\xdef\secsym{\hbox{#1.}}
\bigbreak\bigskip\noindent{\bf Appendix #1. #2}
\writetoca{Appendix {#1.} {#2}}\par\nobreak\medskip\nobreak}
%
%
\def\eqnn#1{\xdef #1{(\secsym\the\meqno)}\writedef{#1\leftbracket#1}%
\global\advance\meqno by1\wrlabel#1}
\def\eqna#1{\xdef #1##1{\hbox{$(\secsym\the\meqno##1)$}}
\writedef{#1\numbersign1\leftbracket#1{\numbersign1}}%
\global\advance\meqno by1\wrlabel{#1$\{\}$}}
\def\eqn#1#2{\xdef #1{(\secsym\the\meqno)}\writedef{#1\leftbracket#1}%
\global\advance\meqno by1$$#2\eqno#1\eqlabel#1$$}
%
%
\def\thm#1{\xdef #1{\secsym\the\mthno}\writedef{#1\leftbracket#1}%
\global\advance\mthno by1\wrlabel#1}
\def\exm#1{\xdef #1{\secsym\the\mexno}\writedef{#1\leftbracket#1}%
\global\advance\mexno by1\wrlabel#1}
%
%
\def\tbl#1{\xdef #1{\secsym\the\tblno}\writedef{#1\leftbracket#1}%
\global\advance\tblno by1\wrlabel#1}
%
\newskip\footskip\footskip14pt plus 1pt minus 1pt 
\def\f@@t{\baselineskip\footskip\bgroup\aftergroup\@foot\let\next}
\setbox\strutbox=\hbox{\vrule height9.5pt depth4.5pt width0pt}
\global\newcount\ftno \global\ftno=0
\def\foot{\global\advance\ftno by1\footnote{$^{\the\ftno}$}}
%
\newwrite\ftfile
\def\footend{\def\foot{\global\advance\ftno by1\chardef\wfile=\ftfile
$^{\the\ftno}$\ifnum\ftno=1\immediate\openout\ftfile=foots.tmp\fi%
\immediate\write\ftfile{\noexpand\smallskip%
\noexpand\item{f\the\ftno:\ }\pctsign}\findarg}%
\def\footatend{\vfill\eject\immediate\closeout\ftfile{\parindent=20pt
\centerline{\bf Footnotes}\nobreak\bigskip\input foots.tmp }}}
\def\footatend{}
%
%
\global\newcount\refno \global\refno=1
\newwrite\rfile
\def\ref{\the\refno\nref}
\def\bref{\nref}
\def\nref#1{\xdef#1{\the\refno}\writedef{#1\leftbracket#1}%
\ifnum\refno=1\immediate\openout\rfile=refs.tmp\fi
\global\advance\refno by1\chardef\wfile=\rfile\immediate
\write\rfile{\noexpand\item{[#1]\ }\reflabel{#1\hskip.31in}\pctsign}\findarg}
\def\findarg#1#{\begingroup\obeylines\newlinechar=`\^^M\pass@rg}
{\obeylines\gdef\pass@rg#1{\writ@line\relax #1^^M\hbox{}^^M}%
\gdef\writ@line#1^^M{\expandafter\toks0\expandafter{\striprel@x #1}%
\edef\next{\the\toks0}\ifx\next\em@rk\let\next=\endgroup\else\ifx\next\empty%
\else\immediate\write\wfile{\the\toks0}\fi\let\next=\writ@line\fi\next\relax}}
\def\striprel@x#1{} \def\em@rk{\hbox{}}
\def\lref{\begingroup\obeylines\lr@f}
\def\lr@f#1#2{\gdef#1{\ref#1{#2}}\endgroup\unskip}

\def\addref#1{\immediate\write\rfile{\noexpand\item{}#1}} 
\def\startrefs#1{\immediate\openout\rfile=refs.tmp\refno=#1}
\def\xref{\expandafter\xr@f}\def\xr@f[#1]{#1}
\def\refs#1{[\r@fs #1{\hbox{}}]}
\def\r@fs#1{\edef\next{#1}\ifx\next\em@rk\def\next{}\else
\ifx\next#1\xref #1\else#1\fi\let\next=\r@fs\fi\next}
%

%
 \newwrite\ffile\global\newcount\figno \global\figno=1
%
%
\def\fig{\the\figno\nfig}
\def\nfig#1{\xdef#1{\secsym\the\figno}%
\writedef{#1\leftbracket \noexpand~\the\figno}%
\ifnum\figno=1\immediate\openout\ffile=figs.tmp\fi\chardef\wfile=\ffile%
\immediate\write\ffile{\noexpand\medskip\noexpand\item{Figure\ \the\figno. }
\reflabel{#1\hskip.55in}\pctsign}\global\advance\figno by1\findarg}
\def\xfig{\expandafter\xf@g}\def\xf@g \penalty\@M\ {}
\def\figs#1{figs.~\f@gs #1{\hbox{}}}
\def\f@gs#1{\edef\next{#1}\ifx\next\em@rk\def\next{}\else
\ifx\next#1\xfig #1\else#1\fi\let\next=\f@gs\fi\next}
%
%
\newwrite\lfile

{\escapechar-1\xdef\pctsign{\string\%}\xdef\leftbracket{\string\{}
\xdef\rightbracket{\string\}}\xdef\numbersign{\string\#}}

\def\writestop{\def\writestoppt{\immediate\write\lfile{\string\pageno%
\the\pageno\string\startrefs\leftbracket\the\refno\rightbracket%
\string\def\string\secsym\leftbracket\secsym\rightbracket%
\string\secno\the\secno\string\meqno\the\meqno}\immediate\closeout\lfile}}
\def\writestoppt{}\def\writedef#1{}

\def\seclab#1{\xdef #1{\the\secno}\writedef{#1\leftbracket#1}\wrlabel{#1=#1}}

\def\subseclab#1{\xdef #1{\secsym\the\subsecno}%
\writedef{#1\leftbracket#1}\wrlabel{#1=#1}}
\def\appsubseclab#1{\xdef #1{\secsym\the\appsubsecno}%
\writedef{#1\leftbracket#1}\wrlabel{#1=#1}}
\def\subsubseclab#1{\xdef #1{\secsym\the\subsecno.\the\subsubsecno}%
\writedef{#1\leftbracket#1}\wrlabel{#1=#1}}
\newwrite\tfile \def\writetoca#1{}
\def\leaderfill{\leaders\hbox to 1em{\hss.\hss}\hfill}
\def\writetoc{\immediate\openout\tfile=toc.tmp
   \def\writetoca##1{{\edef\next{\write\tfile{\noindent ##1
   \string\leaderfill {\noexpand\number\pageno} \par}}\next}}}

%
\catcode`\@=12 
%
%
%
%
%
\def\dbend{{\manual\char127}}
\def\d@nger{\medbreak\begingroup\clubpenalty=10000
    \def\par{\endgraf\endgroup\medbreak} \noindent\hang\hangafter=-2
    \hbox to0pt{\hskip-\hangindent\dbend\hfill}\ninepoint}
\outer\def\danger{\d@nger}

\def\darr#1{\raise1.5ex\hbox{$\leftrightarrow$}\mkern-16.5mu #1}
\def\half{{\textstyle{1\over2}}} 

%
%
\def\al{\alpha}
\def\be{\beta}
\def\ga{\gamma}  
\def\de{\delta}  \def\De{\Delta}
\def\ep{\epsilon}  

\def\et{\eta}
    \def\vth{\vartheta}

\def\la{\lambda} \def\La{\Lambda}
\def\rh{\rho}

  \def\Ps{\Psi}
\def\om{\omega}  
%
%

%

%
%
 
 \def\cD{{\cal D}}
   
\def\cH{{\cal H}}

\def\cW{{\cal W}}

\def\proof{\noindent {\it Proof:}\ }
\def\Box{\hbox{$\rlap{$\sqcup$}\sqcap$}}

%

%
%
\def\amsyes{y }

\def\answ{y }

\ifx\answ\amsyes
\input amssym.def


\def\CC{{\Bbb C}}
\def\ZZ{{\Bbb Z}}
\def\NN{{\Bbb N}}

\def\RR{{\Bbb R}}
\def\bfg{{\frak g}}
\def\bfh{{\frak h}}

\def\hg{{\widehat{\frak g}}}
\def\whg{\hg}
\def\sln{\frak{sl}_N}   
\def\sltw{\frak{sl}_2}  
\def\slth{\frak{sl}_3}

\def\slN{\frak{sl}_N} 

\else
\def\ZZ{{Z\!\!\!Z}}              
\def\CC{{I\!\!\!\!C}}
\def\NN{{I\!\!N}}             

\def\RR{I\!\!R}			
\def\bfg{{\bf g}}
\def\bfh{{\bf h}}

\def\hg{\hat{\bf g}}

\def\sln{s\ell_n}   
\def\sltw{s\ell_2}  
\def\slth{s\ell_3}

\def\slN{s\ell_N}

\fi
%

%
%

%
%

\newsymbol\ltimes 226E
\newsymbol\rtimes 226F
%
%
%

\def\CMP#1{Comm.\ Math.\ Phys.\ {\bf #1}}

\def\LMP#1{Lett.\ Math.\ Phys.\ {\bf #1}}

\def\NPB#1{Nucl.\ Phys.\ {\bf B#1}}

\def\PRep#1{Phys.\ Rep.\ {\bf #1}}

%

%
%
\def\SMu{\hbox{\lower 3pt\hbox{ \epsffile{su10.eps}}}}
\def\SMs{\hbox{\lower 3pt\hbox{ \epsffile{ss10.eps}}}}
\def\SMd{\hbox{\lower 3pt\hbox{ \epsffile{sd10.eps}}}}

\def\SMS{\leavevmode\vadjust{\rlap{\smash%
{\line{{\escapechar=` \hfill\rlap{\hskip.3in%
                 \hbox{\lower 2pt\hbox{\epsffile{sd10.eps}}}}}}}}}}
\def\SMH{\leavevmode\vadjust{\rlap{\smash%
{\line{{\escapechar=` \hfill\rlap{\hskip.3in%
                 \hbox{\lower 2pt\hbox{\epsffile{su10.eps}}}}}}}}}}
%
%
\def\LW#1{\lower .5pt \hbox{$\scriptstyle #1$}}
\def\LWr#1{\lower 1.5pt \hbox{$\scriptstyle #1$}}
\def\LWrr#1{\lower 2pt \hbox{$\scriptstyle #1$}}
\def\RSr#1{\raise 1pt \hbox{$\scriptstyle #1$}}

%


\magnification=1200
\hfuzz=20pt

\hsize=6.5truein
\vsize=9.0truein

\pageno=0
\nopagenumbers


\def\sln{{\frak{sl}}_N}

\def\nroot{k}
\def\sqrtN#1{{\!\!\root{{\nroot}}\of{#1}}}
\def\ni{\noindent}
\def\cox{{\rm h}}
\def\dcox{{\rm h}^\vee}
\def\proof{\noindent {\it Proof.\ }}

\def\bfm{{\bf m}}

\def\eql{~=~}
\def\eqv{~\equiv~}

\def\qeii{QEI$\!$I}
\def\NO#1{~:\!{#1}\!:}

\def\bmu{\bar{\mu}} \def\bep{\bar{\ep}}

\def\cst#1#2{ {\langle #1\rangle \over \langle #2\rangle}}

\font\tencmbsy  = cmbsy10
\newfam\cmbsyfam
\textfont\cmbsyfam=\tencmbsy
\def\bcal#1{{\fam\cmbsyfam\relax#1}}

\font\tencmmib = cmmib10
\font\sevencmmib = cmmib7
\def\bfm#1{\hbox{$#1\textfont1=\tencmmib $}}
\def\sbfm#1{\hbox{$#1\textfont1=\sevencmmib$}}

\def\bfla{\bfm\la}

\def\sbfla{\sbfm\la}



\bref\AFRSa{
J.~Avan, L.~Frappat, M.~Rossi and P.~Sorba,
{\it New ${\cal W}_{q,p}(sl(2))$ algebras from the elliptic algebra
${\cal A}_{q,p}({\hat sl}(2)_c)$}, 
{\tt q-alg/9706013}.}

\bref\AFRSb{
J.~Avan, L.~Frappat, M.~Rossi and P.~Sorba,
{\it Deformed $\cW_N$ algebras from elliptic $sl(N)$ algebras},
{\tt math.QA/9801105}.}

\bref\AKOSa{
H.~Awata, H.~Kubo, S.~Odake and J.~Shiraishi,
{\it Quantum $\cW_N$ algebras and Macdonald polynomials},
\CMP{179} (1996) 401-416,
{\tt q-alg/9508011}.}

\bref\BP{
P.~Bouwknegt and K.~Pilch, 
{\it The deformed Virasoro algebra at roots of unity}, 
{\tt q-alg/9710026}.}

\bref\BS{
P.~Bouwknegt and K.~Schoutens, 
{\it $\cW$ symmetry in conformal field theory},
\PRep{223} (1993) 183-276, {\tt hep-th/9210010}.}

\bref\CPc{
V.~Chari and A.~Pressley,
{\it Quantum affine algebras},
\CMP{142} (1991) 261-283.}

\bref\CPd{
V.~Chari and A.~Pressley,
{\it A guide to quantum groups},
(Cambridge University Press, 1994, Cambridge).}

\bref\DE{
J.~Ding and P.~Etingof,
{\it The center of a quantum affine algebra at the critical level},
Math.\ Res.\ Lett.\ {\bf 1} (1994) 469-480, {\tt hep-th/9403064}.}

\bref\Dr{
V.G.~Drinfel'd, {\it A new realization of Yangians and quantized affine 
algebras}, Sov.\ Math.\ Dokl.\ {\bf 36} (1988) 212-216;
{\it ibid.}, {\it Quantum groups}, proceedings of the ICM, Berkeley, 1986.
}

\bref\FF{
B.~Feigin and E.~Frenkel,
{\it Quantum $\cW$-algebras and elliptic algebras},
\CMP{178} (1996) 653--678, 
{\tt q-alg/9508009}.}

\bref\FJMOP{
B.~Feigin, M.~Jimbo, T.~Miwa, A.~Odesskii and Y.~Pugai,
{\it Algebra of screening operators for the deformed $W_n$ algebra},
{\tt q-alg/9702029}.}

\bref\FRa{
E.~Frenkel and N.~Reshetikhin,
{\it Quantum affine algebras and deformations of the Virasoro 
and $\cW$-algebras},
\CMP{178} (1996) 237--264,
{\tt q-alg/9505025}.}

\bref\FRS{
E.~Frenkel, N.~Reshetikhin and M.A.~Semenov-Tian-Shansky,
{\it Drinfeld-Sokolov reduction for difference operators and
deformations of $\cW$-algebras, I. The case of Virasoro algebra},
{\tt q-alg/9704011}.}

\bref\FRb{
E.~Frenkel and N.~Reshetikhin,
{\it Towards deformed chiral algebras},
{\tt q-alg/9706023}.}

\bref\FRc{
E.~Frenkel and N.~Reshetikhin,
{\it Deformations of $\cal W$-algebras associated to simple Lie algebras},
{\tt q-alg/9708006}.}

\bref\FRaff{
E. Frenkel and N. Reshetikhin,
{\it Deformed $\cW$-algebras and finite-dimensional representations 
of quantum affine algebras}, in preparation.}

\bref\FRe{
E.~Frenkel, talk at the Lie group meeting in Riverside, 2 November 1997,
and private discussions.}

\bref\HYb{
B.-Y.~Hou and W.-L.~Yang,
{\it An $\hbar$-deformation of the $W_N$ algebra and its vertex operators},
J.\ Phys.\ {\bf A}: Math.\ Gen.\ {\bf 30} (1997) 6131-6145,
{\tt hep-th/9701101}.}

\bref\JLMP{
M.~Jimbo, M.~Lashkevich, T.~Miwa and Y.~Pugai,
{\it Lukyanov's screening operators for the deformed Virasoro algebra},
Phys.\ Lett.\ {\bf 229A} (1997) 285-292,
{\tt hep-th/9607177}.}

\bref\JS{
M.~Jimbo and J.~Shiraishi,
{\it A coset-type construction for the deformed Virasoro algebra},
{\tt q-alg/9709037}.}

\bref\Kog{
A.~Kogan,
{\it Deformations of the classical W-algebras associated to 
$D_n$, $E_6$ and $G_2$}, {\tt q-alg/9710027}.}

\bref\LPa{
S.~Lukyanov and Y.~Pugai,
{\it Bosonization of ZF algebras: Direction toward deformed Virasoro
algebra}, J.\ Exp.\ Theor.\ Phys.\ {\bf 82 } (1996) 1021-1045,
{\tt hep-th/9412128}.}

\bref\LPb{
S.~Lukyanov and Y.~Pugai,
{\it Multi-point local height probabilities in the integrable RSOS
model}, \NPB{473} (1996) 631-658,
{\tt hep-th/9602074}.}

\bref\MD{
I.G.~Macdonald,
{\it Symmetric functions and Hall polynomials},
(Oxford University Press, Oxford, 1995).}

\bref\RSTS{
N.~Reshetikhin and M.~Semenov-Tian-Shansky,
{\it Central extensions of quantum current groups},
\LMP{19} (1990) 133-142.}

\bref\SS{
A.M.~Semenov-Tian-Shansky and A.V.~Sevostyanov,
{\it Drinfeld-Sokolov reduction for difference operators and
deformations of $\cW$-algebras, II. General semisimple case},
{\tt q-alg/9702016}.}

\bref\SKAO{
J.~Shiraishi, H.~Kubo, H.~Awata and S.~Odake,
{\it A quantum deformation of the Virasoro algebra and the Macdonald
symmetric functions},
\LMP{38} (1996) 33--51,
{\tt q-alg/9507034}.}

%
%
\line{}
\vskip2cm
\centerline{\bf ON  DEFORMED $\bcal{W}$-ALGEBRAS}\smallskip
\centerline{\bf AND QUANTUM AFFINE ALGEBRAS}\bigskip
\vskip1cm

\centerline{Peter BOUWKNEGT$\,^{1}$ and Krzysztof PILCH$\,^{2}$}
\bigskip

\centerline{\sl $^1$ Department of Physics and Mathematical Physics}
\centerline{\sl University of Adelaide}
\centerline{\sl Adelaide, SA~5005, AUSTRALIA}
\bigskip

\centerline{\sl $^2$ Department of Physics and Astronomy}
\centerline{\sl University of Southern California}
\centerline{\sl Los Angeles, CA~90089-0484, USA}
\medskip
\vskip1.5cm

\centerline{\bf ABSTRACT}\bigskip

\vbox{\leftskip 2.0truecm \rightskip 2.0truecm
\noindent 
We discuss some aspects of the deformed $\cW$-algebras
$\cW_{q,t}[\bfg]$.  In particular, we derive an explicit formula for
the Kac determinant, and discuss the center when $t^2$ is a primitive
$k$-th root of unity.  The relation of the structure of
$\cW_{q,t}[\bfg]$ to the representation ring of the quantum affine
algebra $U_q(\widehat \bfg)$, as discovered recently by Frenkel and
Reshetikhin, is further elucidated in some examples. }

\vfil
\leftline{ADP-97-49/M60}
\leftline{USC-97/021}
\line{{\tt math.QA/9801112}\hfil January 1998}

\eject


\baselineskip=1.4\baselineskip

\footline{\hss \tenrm -- \folio\ -- \hss}

\newsec{Introduction}

In recent years there has been a considerable interest in
understanding the role of infinite dimensional quantum algebras in the
theory of off-critical integrable models of statistical mechanics. In
particular, it has become clear that the algebraic framework of those
theories closely parallels that of their critical counterpart -- the
conformal field theories.

In this context the algebras of
particular interest are the so-called deformed $\cW$-algebras
$\cW_{q,t}[\bfg]$ associated to arbitrary simple finite dimensional 
Lie algebras $\bfg$.  They can be considered to be deformations of the 
$\cW$-algebras $\cW[\bfg]$ that arise in conformal field theory.%
\foot{For a review of $\cW$-algebras $\cW[\bfg]$ see, e.g., [\BS] and
the references therein.}
In analogy with the 
undeformed case, the algebras $\cW_{q,t}[\bfg]$ 
are most conveniently defined as the centralizer of a set of screening
operators $S_i^\pm$ in a deformed Heisenberg algebra $\cH_{q,t}[\bfg]$.

The simplest example, the deformed Virasoro algebra ${\rm Vir}_{q,t}\cong
\cW_{q,t}[\sltw]$, was introduced in [\FRa,\SKAO] and further studied in
[\FRS,\SS,\LPb,\AFRSa,\JLMP,\JS,\BP].  In [\LPa,\LPb] it has been argued  
that this deformed Virasoro algebra plays the role of the dynamical symmetry
algebra in the Andrews-Baxter-Forrester RSOS models.  The higher rank
generalizations, $\cW_{q,t}[\sln]$, were introduced in 
[\FF,\AKOSa,\HYb,\FJMOP,\FRc,\AFRSb]. 
The deformed $\cW$-algebras $\cW_{q,t}[\bfg]$, for arbitrary simple finite
  dimensional Lie algebras $\bfg$, were introduced recently by Frenkel and
  Reshetikhin [\FRc] and further studied in, e.g., [\Kog].

In this paper we further investigate the structure of $\cW_{q,t}[\bfg]$.
In particular, we derive explicit expressions for elements in the 
centralizer of the screening operators for $\bfg=\sltw$, i.e., generators
of $\cW_{q,t}[\sltw]$, and formulate an algorithm to obtain generators
of $\cW_{q,t}[\bfg]$ by `pasting' of the various $\sltw$ directions.  
We illustrate this procedure in the case of the rank $2$ simple Lie algebras.
This explicit construction further elucidates the connection of 
$\cW_{q,t}[\bfg]$ to the Grothendieck ring ${\rm Rep\,}(U_q(\whg))$ of
finite dimensional representations of the quantum affine algebra $U_q(\whg)$
as discovered in [\FRc].  

The main ideas and some of the formulae 
of this part of the paper are already, at least implicitly, contained in
[\FRc].  The main results of this paper are a formula for the determinant 
of the contravariant bilinear form (`Shapovalov form') on the Verma modules
over $\cW_{q,t}[\bfg]$ (Theorem 4.4) and an explicit
construction of the center of $\cW_{q,t}[\bfg]$
when one of the deformation parameters $t^2$ is a primitive $\nroot$-th
root of unity (Theorem 5.1).  
These two results generalize our results for $\bfg=\sltw$ in [\BP].

The paper is organized as follows: In Section 2 we first recall,
following [\FRc], the definition of the deformed $\cW$-algebra
$\cW_{q,t}[\bfg]$ and formulate a conjecture for the structure of its
generators $T_i(z)$.  This conjecture (Assumption 2.4) 
has been verified in a number of
cases, including $\bfg=\sln$ and all rank $2$ simple Lie algebras.  
It is a slight refinement of Conjecture 1 in [\FRc] and will
be the starting point for the analysis in this paper.
We also prove some simple corollaries of this conjecture which will be needed
in the discussion of the main results of this paper.
In Section 3 we provide further evidence for the conjectured structure 
of the $\cW_{q,t}[\bfg]$ generators by explicitly working out the centralizer
of the screening operators in the case of $\bfg=\sltw$.  We also formulate
an algorithm to obtain generators
of $\cW_{q,t}[\bfg]$ by `pasting' of the various $\sltw$ directions, and
further elucidate the connection of those results to the representation
theory of quantum affine algebras.  In Section 4 we prove a formula for
the Kac determinant of $\cW_{q,t}[\bfg]$ and in Section 5 we discuss the
center of $\cW_{q,t}[\bfg]$ for $t^2$ a primitive $\nroot$-th root of unity.
In an appendix we illustrate several of the issues raised above in the 
case of the deformed algebras $\cW_{q,t}[\bfg]$ where $\bfg$ is a 
simple finite dimensional Lie algebra of rank $2$.  In particular, we
give explicit formulae for all the generators $T_i(z)$ and their commutation
relations.

\newsec{Deformed $\bcal{W}$-algebras }

In this section we review the construction of the deformed algebra
$\cW_{q,t}[\bfg]$, where $\bfg$ is an arbitrary simple Lie algebra of
rank $\ell$, and $q,t\in\CC$ are deformation parameters. The first
part of this section closely follows [\FRc], which the reader should
consult for further details.
\thm\thSTWH
\proclaim Definition \thSTWH\ [\FRc]. The deformed Heisenberg algebra 
$\cH_{q,t}[\bfg]$ is the associative algebra with the (root type)
generators $a_i[n]$, $i=1,\ldots,\ell$, $n\in\ZZ$, satisfying the
relations
\eqn\STWaa{
[a_i[m],a_j[n]]\eql {1\over m}(q^m-q^{-m})(t^m-t^{-m})
B_{ij}(q^m,t^m)\de_{m+n,0}\,,
}
where $B(q,t)=D(q,t)C(q,t)$ and%
\foot{\rm We use the
standard $q$-notation, $[n]_q=(q^n-q^{-n})/(q-q^{-1})$.} 
\eqn\STWac{
C_{ij}(q,t)\eql(q^{r_i}t^{-1}+q^{-r_i}t)\de_{ij}-[I_{ij}]_q\,,
}
\eqn\STWad{
D_{ij}(q,t)\eql [r_i]_q\de_{ij}\,,
}
are, respectively,  deformations of the Cartan matrix $C=(C_{ij})$
and the diagonal matrix $D={\rm diag}(r_1,\ldots,r_n)$ of $\bfg$.
\medskip

We recall that $I_{ij}=2\de_{ij}-C_{ij}$ is the incidence matrix of
$\bfg$, while the 
relatively prime integers $r_i$ are given in terms of lengths of
canonically normalized simple roots by $r_i=r^\vee (\al_i,\al_i)/2$,
where $r^\vee$ is the dual tier number of $\bfg$, i.e., the maximal number
of edges connecting two vertices of the Dynkin diagram. In particular,
when $\bfg$ is simply laced, we have $r^\vee=1$, $r_i=1$, and
\eqn\STWae{
B_{ij}(q,t)\eql C_{ij}(q,t)\eql (qt^{-1}+q^{-1}t)\de_{ij}-I_{ij}\,.
}

The set of fundamental weight type generators, $y_i[n]$,
$i=1,\ldots,\ell$, $n\in\ZZ$, is defined by 
\eqn\STWaf{
a_i[m]\eql\sum_j C_{ji}(q^m,t^m)y_j[m]\,.}
They satisfy
\eqn\STWaga{
[a_i[m],y_j[n]] \eql {1\over m} (q^{r_im}-q^{-r_im})(t^m-t^{-m}) \de_{ij}
  \de_{m+n,0}\,,
}
and
\eqn\STWag{
[y_i[m],y_j[n]]\eql {1\over m}
(q^m-q^{-m})(t^m-t^{-m}) M_{ij}(q^m,t^m)\de_{m+n,0}\,,
}
where $M(q,t)=D(q,t)C(q,t)^{-1}$.  An explicit formula for the matrices 
$M(q,t)$ for Lie algebras of classical type can be found in Appendix C
of [\FRc].

In the following we will use the generating series
\eqn\STWah{
A_i(z)\eql t^{2(\rh^\vee,\al_i)}q^{-2r^\vee(\rh,\al_i)+2a_i[0]}
\NO{ \exp \Bigg(\sum_{m\not=0} a_i[m]z^{-m}\Bigg) }\,,}
and
\eqn\STWai{
Y_i(z)\eql t^{2(\rh^\vee,\om_i)}q^{-2r^\vee(\rh,\om_i)+2y_i[0]}
\NO{ \exp \Bigg(\sum_{m\not=0} y_i[m]z^{-m}\Bigg)}\,,}
where $\om_i$ are the fundamental weights and $\rh$ is 
the Weyl vector of $\bfg$, i.e.,
$(\om_i,\al_j^\vee) = \de_{ij}$ and $(\rh,\al_i^\vee) = 1$. 

A Fock module, $F(\mu)$, of $\cH_{q,t}[\bfg]$, where $\mu$ is a weight
of $\bfg$, is freely generated by $a_i[m]$, $m<0$, from the vacuum
state $|\mu\rangle$, which satisfies $a_i[0]|\mu\rangle=
(\mu,\al_i)|\mu\rangle$ and $a_i[m]|\mu\rangle=0$, $m>0$.  
It decomposes as
$F(\mu)\eql \coprod_{n\geq 0} F(\mu)_{(n)}
$
under the action of, $d$, the derivation of $\cH_{q,t}[\bfg]$ defined by
\eqn\STTxb{
[d,a_i[m]]\eql-ma_i[m]\,,\quad m\in\ZZ\,.
} 
There are two canonical bases of 
$F(\mu)_{(n)}$ consisting of vectors
\eqn\STTxc{
o[-\bfla]|\mu\rangle
~\equiv~o_1[-\la^{(1)}_1]\ldots o_1[-\la^{(1)}_{l_1}]\ldots
o_\ell[-\la^{(\ell)}_1]\ldots o_\ell[-\la^{(\ell)}_{l_\ell}]|\mu\rangle\,,
}
where $o_i[m]=a_i[m]$ and $o_i[m]=y_i[m]$, respectively, and
$\bfla\vdash n$ runs over all multi-partitions of $n$, i.e.\
$\bfla=(\la^{(1)},\ldots,\la^{(\ell)})$ with $|\bfla|=\sum_i |\la^{(i)}|=n$.

The formal power series generated by terms of the form
$$
\NO{ \partial_z^{n_1}Y_{i_1}(zq^{a_1}t^{b_1})^{\ep_1}\ldots 
\partial_z^{n_l}Y_{i_l}(zq^{a_l}t^{b_l})^{\ep_l}}\,,
$$
where $\ep_i=\pm1$, $n_i\geq 0$ and $(a_i,b_i)\in L\subset \ZZ\times
\ZZ$, together with the Fock module $F(0)$, form a deformed chiral
algebra  in the sense of [\FRb]. It will be denoted by
${\bf H}_{q,t}[\bfg]$.

To construct maps between Fock modules we extend $\cH_{q,t}[\bfg]$ to
$\cH'_{q,t}[\bfg]$ by introducing operators $Q_i$, $i=1,\ldots,\ell$,
such that $e^{Q_i}$ are shift operators satisfying
\eqn\STWaj{
[a_i[m],e^{Q_j}]\eql B_{ij}\be \de_{m,0}e^{Q_j}\,,
}
or, equivalently,
\eqn\STWaja{
[y_i[m], e^{Q_j}]\eql r_i \be \de_{ij} \de_{m,0} e^{Q_j}\,,
}
where $\be$ is given by $t=q^\be$ and $B_{ij}=B_{ij}(1,1)$. Now, the
screening currents, $S_i^\pm(z)$, are defined as the generating series
\eqn\STWak{
S_i^+(z)\eql e^{-Q_i/r_i}z^{-s_i^+[0]} \NO{ \exp \Bigg(\sum_{m\not=0}
s_i^+[m]z^{-m}\Bigg) } \,,
}
\eqn\STWal{
S_i^-(z)\eql e^{Q_i/\be}z^{s_i^-[0]} \NO{ \exp \Bigg(-\sum_{m\not=0}
s_i^-[m]z^{-m}\Bigg) }\,,
}
where
\eqn\STWam{
s_i^+[m]\eql {a_i[m]\over q^{mr_i}-q^{-mr_i}}\,,\quad m\not=0,\qquad
s_i^+[0]\eql {a_i[0]\over r_i}\,,}
\eqn\STWam{
s_i^-[m]\eql {a_i[m]\over t^{m}-t^{-m}}\,,\quad m\not=0,\qquad
s_i^-[0]\eql {a_i[0]\over \be}\,.}
The screening operators $S_i^+:F(0)\rightarrow F(-\be \al_i^\vee)$
and $S^-_i:F(0)\rightarrow F(r^\vee\al_i)$ are
defined by $S_i^{\pm}=S_i^\pm[1]$, where $S_i^\pm(z) = \sum_{m\in\ZZ}
S_i^{\pm}[m] z^{-m}$.
\thm\thSTab
\proclaim Definition \thSTab\ [\FRc]. The deformed $\cW$-algebra 
$\cW_{q,t}[\bfg]$ is the associative algebra, topologically generated
by the Fourier coefficients of the maximal deformed chiral subalgebra
of\/ ${\bf H}_{q,t}[\bfg]$, which commutes with the screening operators
$S_i^\pm$, $i=1,\ldots,\ell$.
\medskip

An explicit construction of the generators of $\cW_{q,t}[\bfg]$ has
been carried out completely in the case $\bfg=\slN$ [\FF,\AKOSa] and,
partially, when $\bfg$ is a Lie algebra of classical type [\FRc],
where, in particular, the generator
$T_1(z)$ was constructed for all classical Lie algebras $g$.
Additional results are known for the classical limit $t\to1$, i.e.\
for the Poisson algebra $\cW_{q,1}[\bfg]$ (see, in particular,
[\FRa,\FRc,\Kog]).
In Appendix A we give a complete result when $\bfg$ is one of the
rank $2$ simple Lie algebras.  

To elucidate the structure of those generators, let us first consider
the case $\bfg=\slN$.
\thm\thSTac
\proclaim Theorem \thSTac\ [\SKAO,\FF,\AKOSa,\FRc]. The 
algebra $\cW_{q,t}[\sln]$ is generated by the Fourier modes of
fields 
\eqn\STWbb{
T_i(z)\eql\sum_{1\leq j_1<\ldots<j_i\leq N} \NO{ \La_{j_1}(zp^{-i+1})
\La_{j_2}(zp^{-i+3})\ldots 
\La_{j_i}(zp^{i-1})}\,,\qquad i=1,\ldots,N-1\,,}
where $\La_i(z)$ are defined recursively by
\eqn\STWbc{\eqalign{
\La_1(z)& \eql Y_1(z)\,,\cr
\La_i(z)& \eql \NO{ \La_{i-1}(z)A_{i-1}(zp^{-i+1})^{-1}}\,,\qquad
i=2,\ldots,N\,,\cr}
}
and $p=qt^{-1}$. The $\La_i(z)$ satisfy 
\eqn\STWyz{
\NO{ \La_1(z) \La_2(zp^2) \ldots \La_N(zp^{2(N-1)}) } \eql 1\,.
}
\medskip

In the remaining cases the structure of the generators of
$\cW_{q,t}[\bfg]$ can be motivated by the explicit examples, various
limiting cases (such as the conformal limit $q\rightarrow 1$,
$\be={\rm const}$) in which the algebra is known, and, most interestingly, a
natural interpretation of formulae like \STWbb\ as characters of
finite dimensional irreducible representations of 
the quantum affine algebra $U_q(\widehat\bfg)$, where $\widehat\bfg$ is the
affine Lie algebra corresponding to $\bfg$. 
A conjecture was first formulated in [\FRc]; in order to prove 
the main results of this paper we will need a slightly sharper version
of that conjecture.  Henceforth, we will assume that the 
following holds:
\thm\thSTad
\proclaim Assumption \thSTad. The algebra $\cW_{q,t}[\bfg]$
is generated by the Fourier modes $T_i[m]$ 
of the fields $T_i(z)=\sum_{m\in\ZZ} T_i[m]z^{-m}$, $i=1,\ldots,\ell$,
\eqn\STWbe{
T_i(z)\eql \sum_{{\la\in P(V(\om_i))}} 
  \sum_{k_\la=1}^{{\rm mult\,}\la} \ c^{\om_i,(k_\la)}_{\la}(q,t) \,
  Y^{\om_i,(k_\la)}_\la(z)\,,
}
where $\la$ runs over the weights $P(V(\om_i))$ of
the finite dimensional irreducible representation $V(\om_i)$ of
$U_q(\widehat\bfg)$ with highest weight $\om_i$, and each 
$Y^{\om_i,(k_\la)}_\la(z)$, with $\la=\om_i-\sum \al_{i_j}$, is of the form
\eqn\STWbea{
Y^{\om_i,(k_\la)}_\la(z) \eql \NO{ Y_i(z)A_{i_1}(zq^{a_1}t^{b_1})^{-1}
\ldots A_{i_k}(zq^{a_k}t^{b_k})^{-1} }\,,
}
for some choice of integers $a_i,b_i\in\ZZ$.
Furthermore, we assume that \STWbe\ is obtained by the algorithm of pasting  
$\sltw$ directions (see Section 3). In particular, $T_i(z)$ 
is uniquely determined by its highest weight component $Y^{\om_i}_{\om_i}(z)
=Y_i(z)$.
We normalize $T_i(z)$ by choosing $c^{\om_i}_{\om_i}(q,t)=1$.
\medskip

\noindent {\bf Remark.}  
For $V(\om_1)$ the set of weights is the same as for the 
irreducible $\bfg$ module of highest weight $\om_1$. However, in
general, this set is bigger, except for $\sln$ (see, e.g., [\CPd]
for background material on quantum affine algebras). \medskip

Again, we emphasize that Assumption \thSTad\ holds in all known cases
including, in particular, $\bfg = \sln$ [\FF,\AKOSa,\FRc] and all
rank $2$ simple Lie algebras $\bfg$ (Appendix A).
In Section 3 we will generalize it to arbitrary finite dimensional
irreducible $U_q(\whg)$ modules $V$.
\thm\thSTTfx
\proclaim Theorem \thSTTfx. The coefficients
$c^{\om_i,(k_\la)}_\la(q,t)$ defined by \STWbe\  
satisfy
\eqn\STTyg{
  c^{\om_i,(k_\la)}_\la(q^{-1},t^{-1}) \eql c^{\om_i,(k_\la)}_\la(q,t)\,,
}
and
\eqn\STTxk{
\sum_{k_\la=1}^{{\rm mult}\,\la}\ c^{\om_i,(k_\la)}_{w\la}(q,t) \eql 
  \sum_{k_\la=1}^{{\rm mult}\,\la}\ 
  c^{\om_i,(k_\la)}_\la(q,t) \,,\qquad w\in W\,.
} 
In particular, if ${\rm mult\,}\la=1$, then the coefficients 
$c_\la^{\om_i}(q,t)$ are invariant with respect to the action of the Weyl
group $W$ of $\bfg$.  In addition, for simply-laced $\bfg$, we have 
\eqn\STTyga{
c^{\om_i,(k_\la)}_\la(t,q) \eql c^{\om_i,(k_\la)}_\la(q,t)\,.
}
\medskip

Since the proof requires a more detailed analysis of the structure of
$T_i(z)$, we defer it to Section 3.

Note that the matrices $C_{ij}(q,t)$, $B_{ij}(q,t)$ and $M_{ij}(q,t)$
are all invariant under $(q,t) \to (q^{-1},t^{-1})$.  This implies that
$\cH_{q,t}[\bfg]\cong\cH_{q^{-1},t^{-1}}[\bfg]$.  In fact, let us combine 
this transformation with a $\ZZ_2$-automorphism of $\cH_{q,t}[\bfg]$ and define
\eqn\STWbf{ \eqalign{
\vth(q) & \eql q^{-1} \,,\qquad \vth(t) \eql t^{-1} \,,\cr
\vth(q^{a_i[0]}) & \eql q^{-a_i[0]}\,,\qquad
\vth(a_i[m]) \eql -a_i[m]\,,\quad m\neq0\,,\cr}
}
then 
\eqn\STWbg{
\vth(A_i(z)) \eql A_i(z)^{-1}\,,\qquad \vth(Y_i(z)) \eql Y_i(z)^{-1}\,,
} 
while the screening currents \STWak\ and \STWal\ are invariant under 
$\vth$ provided we define 
\eqn\STWbh{
\vth( Q_i ) \eql Q_i \,.
}
The invariance of the screening operators implies that $\cW_{q,t}[\bfg]$
is invariant under $\vth$.  Indeed, in Section 3 we prove
\thm\thSTae
\proclaim Theorem \thSTae\ [\FRc]. Under the transformation $\vth$, 
\eqn\STWbi{
\vth( T_i(z) ) \eql T_{i^*}(zq^{r^\vee\dcox}t^{-\cox})\,,
}
where $T_{i^*}(z)$ is the generator corresponding to the weight 
$\om_i^* = -w_0 \om_i$, conjugate to $\om_i$.
\medskip

The fields $T_i(z)$ satisfy the exchange relations [\FF,\FRc]
\eqn\STWdaa{
T_i(z)T_j(w)\eql S_{T_i,T_j}({w\over z}) T_j(w)T_i(z)\,,}
where $S_{T_i,T_j}({x})\eql f_{ij}(x)^{-1}f_{ij}({1/ x})$ and  
$f_{ij}(x)$ is defined by
\eqn\STFbx{
Y_i(z)Y_j(w)\eql f_{ij}({w\over z})^{-1}\NO{ Y_i(z)Y_j(w) }\,.
}
A straightforward calculation using \STWag\ yields
\eqn\STWda{
f_{ij}(x)\eql\exp \Bigg(-\sum_{m>0}
(q^m-q^{-m})(t^m-t^{-m})M_{ij}(q^m,t^m){x^m\over m}\Bigg)\,.}
{}Finally,
the products of operators in \STWda\ are understood in the sense of
analytic continuation.
Using standard techniques one can derive from
\STWda\ the corresponding commutation relations for the modes $T_i[m]$
(see [\FF,\AKOSa] for the case of $\sln$ and Appendix A for the 
rank $2$ simple Lie algebras).

The Verma module, $M(h)$, of $\cW_{q,t}[\bfg]$ is defined as usual
[\FF,\AKOSa]. It is generated from the highest weight state,
$|h\rangle$, satisfying $T_i[0]|h\rangle=h_i|h\rangle$ and
$T_i[m]|h\rangle=0$, $m>0$, and  decomposes under the action of $d$ as
$M(h)=\coprod_{n\geq0}M(h)_{(n)}$. A basis of $M(h)_{(n)}$ consists
of vectors
\eqn\STTxh{
T[-\bfla] |\mu\rangle ~\equiv~T_1[-\la^{(1)}_1]\ldots
T_1[-\la^{(1)}_{l_1}]\ldots T_\ell[-\la^{(\ell)}_1]\ldots
T_\ell[-\la^{(\ell)}_{l_\ell}]|\mu\rangle\,, } 
indexed by multi-partitions $\bfla$. 

The realization of $T_i(z)$ in terms of $y_i(z)$ defines a
homomorphism $\imath:M(h(\mu)) \rightarrow F(\mu)$, uniquely
determined by $\imath(|h(\mu)\rangle)=|\mu\rangle$. The highest weight,
$h(\mu)$, of the Verma module can be found by evaluating $T_i[0]$ on
the highest weight vector of $F(\mu)$,
$T_i[0]|\mu\rangle=h_i(\mu)|\mu\rangle$.  As a consequence of 
Assumption \thSTad\ we have
\thm\thTTWxa
\proclaim Corollary \thTTWxa.  The eigenvalues $h_i(\mu)$  are given by
\eqn\STTxj{
h_i(\mu) \eql \sum_{\la\in P(V(\om_i))}\ \biggl(\ 
  \sum_{k_\la =1}^{{\rm mult\,}\la} \ c^{\om_i,(k_\la)}_\la(q,t) \biggr) \, 
  q^{-2r^\vee (\rh,\la)} t^{2(\rh^\vee,\la)} q^{2(\mu,\la)}\,,
}
where the sum runs over the weights  of the irreducible module 
$V(\om_i)$ of $U_q(\widehat{\bfg})$.
\medskip

It follows from Theorem \thSTTfx\ that the $h_i(\mu)$ are invariant
under 
\eqn\TSSaa{
\mu ~\rightarrow~ w*\mu ~\equiv~ w(\mu - r^\vee \rh +\be \rh^\vee) + 
  (r^\vee \rh -\be \rh^\vee)\,,
}
for each $w$ in the Weyl group $W$ of $\bfg$.  In fact,
\thm\thTSaa
\proclaim Lemma \thTSaa.  For generic $q,t\in\CC$, we have $h_i(\mu) = 
h_i(\mu')$ for all $i=1,\ldots,\ell$, if and only if there exists 
a $w\in W$ such that $\mu' = w*\mu$.\medskip

For future use, also note that
\eqn\TSSab{
h_i(\bmu) \eql h_{i^*}(\mu)\,,
}
where 
\eqn\TSSac{
\bmu \eqv -\mu+2r^\vee\rh-2\be\rh^\vee\,.
}


\newsec{Explicit generators of $\cW_{q,t}[\bfg]$}

In this section we outline an algorithm to explicitly compute elements 
in the commutant of the screening operators $S_i^\pm$, i.e., generators
of $\cW_{q,t}[\bfg]$.
First we analyze the case $\bfg=\sltw$ in detail and then 
indicate how to obtain the result for arbitrary $\bfg$ by pasting together
the various $\sltw$ directions.  
In Appendix A we illustrate this procedure in
the case of the fundamental generators of $\cW_{q,t}[\bfg]$, for all 
rank $2$ simple Lie algebras $\bfg$.
The main ideas of this section, as well as some of the explicit formulae,
are already contained in [\FRb,\FRc] -- our main motivation is to make 
these ideas more explicit in order to provide further support for 
Assumption \thSTad, and to prove Theorems \thSTTfx\ and \thSTae\
which will play a crucial role in our analysis of the Kac determinant 
for $\cW_{q,t}[\bfg]$.

{}From the commutation relations \STWaga\ and \STWaja\ 
it follows%
\foot{The normal ordering includes the prescription 
to put the $y_i[0]$ to the right of the $Q_j$.}
\eqn\eqAPa{ \eqalign{
Y_i(z) S_i^+(w) & \eql t^{-2} \left( { 1-  t {w\over z}\over 
  1-  t^{-1} {w\over z}} \right) \NO{ Y_i(z) S_i^+(w)  } \,,\cr
Y_i(z) S_i^-(w) & \eql q^{2r_i} \left( { 1- q^{-r_i} {w\over z}  \over 
  1-  q^{r_i} {w\over z} } \right) \NO{ Y_i(z) S_i^+(w)  } \,,\cr
S_i^+(w) Y_i(z) & \eql \left( { 1- t^{-1}  {z\over w}\over 
  1- t {z\over w} } \right) \NO{ Y_i(z) S_i^+(w)  } \,,\cr
S_i^-(w) Y_i(z) & \eql \left( { 1- q^{r_i} {z\over w}\over 
  1- q^{-r_i}  {z\over w}} \right) \NO{ Y_i(z) S_i^-(w)  } \,,\cr}
}
while, for $i\neq j$,
\eqn\eqAPb{ \eqalign{
Y_i(z) S_j^\pm(w) & \eql \NO{ Y_i(z) S_j^\pm(w)  } \,,\cr
S_j^\pm(w) Y_i(z) & \eql \NO{ Y_i(z) S_j^\pm(w)  } \,.\cr}
}
Furthermore, we have the difference relations [\FRc]
\eqn\eqAPc{ \eqalign{
S_i^+(z) & \eql q^{2r_i} t^{-2} \NO{ A_i(zq^{r_i}) S_i^+(zq^{2r_i}) }\,,\cr
S_i^-(z) & \eql q^{2r_i} t^{-2} \NO{ A_i(zt^{-1}) S_i^-(zt^{-2}) }\,.\cr}
}
Our aim is to find certain combinations of vertex operators such that,
by using the difference relations \eqAPc, the commutator
with the screening currents $S_i^+(w)$ and 
$S_i^-(w)$ can be written as total $q^{2r_i}$- and $t^2$-differences, 
respectively.
We recall the definition of a total $a$-difference $\cD_a$,
\eqn\eqAPd{
\cD_a\cdot f(w) \eql { f(w) - f(wa) \over w(1-a)}\,.
}
This then obviously implies that those combinations are in the 
commutant of $S_i^\pm = S_i^\pm[1]$.

Now consider $\bfg=\frak{sl}_2$, and a vertex operator of the form
(we write $Y(z)=Y_1(z)$ and $A(z)=A_1(z)$ in the case of $\frak{sl}_2$)
\eqn\eqAPe{
\Ps(z) \eql \NO{ Y(z_1) \ldots Y(z_m)}\,.
}
For our purposes it will be sufficient to consider the case 
\eqn\eqAPf{
z_i \eql zq^{2a_i} t^{2b_i}\,,\qquad i=1,\ldots,m\,,
}
for some $({\bf a},{\bf b}) = ((a_1,b_1),\ldots,(a_m,b_m))\in\RR^m
\times\RR^m$.
For convenience we define 
\eqn\eqAPg{
\xi_{ij} \eql \sqrt{{z_i\over z_j}} \eql q^{a_i-a_j}t^{b_i-b_j}
  \eql  \xi_{ji}^{-1}\,.
}
We first consider the case when the $z_i$ are in generic position, i.e.,
$\xi_{ij} \not\in q^{\ZZ} t^{\ZZ}$ for all pairs $(i,j)$, and take the 
non-generic limit of the resulting expressions at a later stage.
Let us introduce
\eqn\eqAPh{
Y^{(+)}(z) \eql Y(z)\,,\qquad Y^{(-)}(z) \eql \NO{Y(z) A(zq^{-1}t)^{-1}} \eql 
  Y(zq^{-2} t^2)^{-1}\,,
}
and consider the $2^m$ vertex operators
\eqn\eqAPh{
\Ps_{(\ep_1\ldots\ep_m)}(z) \eql \NO{ Y^{(\ep_1)}(z_1) \ldots 
  Y^{(\ep_m)}(z_m)}\,,
}
in particular $\Ps(z)=\Ps_{(++\ldots+)}(z)$.  A tedious, but straightforward,
calculation yields
\eqn\eqAPi{ \eqalign{
[& \Ps_{(\ep_1\ldots\ep_m)}(z)  ,S^-(w)] \cr
& \eql q^{m_+ - m_-}
  (q-q^{-1})\sum_i {\rm sgn}({\ep_i})
  \Biggl( \prod_{j\neq i} \et_{ij}^{(\ep_i\ep_j)}(q,t) \Biggr)
  \de( x_i^{(\ep_i)} q) \NO{\Ps_{(\ep_1\ldots\ep_m)}(z) S^-(w) }\,,\cr}
}
where 
\eqn\eqAPj{
M_\pm \eql \{ i\ |\ \ep_i =\pm \}\,, \qquad m_\pm \eql |M_\pm|\,,
}
\eqn\eqAPk{
x_i^{(+)} \eql x_i \eql {w\over z_i}\,,\qquad x_i^{(-)} \eql x_i t^{-2}  \,,
}
and
\eqn\eqAPla{
\et_{ij}^{(\ep_i,\ep_j)}(q,t) \eqv \et^{(\ep_i,\ep_j)}(q,t;\xi_{ij})\,,
}
with
\eqn\eqAPl{ \eqalign{
\et^{(++)}(q,t;\xi) & \eql {q \xi^{-1} -q^{-1} \xi \over \xi^{-1}
  -\xi } \,,\cr
\et^{(-+)}(q,t;\xi) & \eql {q t^{-1} \xi^{-1} - q^{-1}t \xi \over 
  t^{-1}\xi^{-1}  - t\xi } \,,\cr
\et^{(+-)}(q,t;\xi) & \eql {q t^{-1} \xi - q^{-1}t \xi^{-1} \over 
  t^{-1}\xi  - t\xi^{-1} } \,,\cr
\et^{(--)}(q,t;\xi) & \eql {q \xi -q^{-1} \xi^{-1} \over \xi
  -\xi^{-1} } \,.\cr}
}
The coefficients $\et^{(\ep,\ep')}(q,t;\xi)$ have the following 
easily verifiable properties
\eqn\eqAPm{ \eqalign{
\et^{(\ep,\ep')}(q,t;\xi^{-1}) & \eql \et^{(\bep,\bep')}(q,t;\xi)\,, \cr
\et^{(\ep,\ep')}(q^{-1},t^{-1};\xi) & \eql \et^{(\bep,\bep')}(q,t;\xi) \,,\cr
\et^{(\ep,\ep')}(q,t;(qt^{-1})^{\ep'}\xi) & \eql 
  \et^{(\ep,\bep')}(q,t;\xi)^{-1} \,,\cr}
}
where $\bep=\mp$ for $\ep=\pm$.
\thm\thAPa
\proclaim Theorem \thAPa.  For parameters $z_i=zq^{2a_i}t^{2b_i}$,
$({\bf a},{\bf b}) = ((a_1,b_1),\ldots,((a_m,b_m)) \in \RR^m \times \RR^m$,
in generic position (i.e., $\xi_{ij}\not\in q^{\ZZ}t^{\ZZ}$), the 
vertex operators 
\eqn\eqAPn{
T_{({\bf a},{\bf b})}(z) \eql \sum_{\ep_i} \ 
  \ga_{(\ep_1\ldots\ep_m)}(q,t)\, \Ps_{(\ep_1\ldots\ep_m)}(z) \,,
}
where
\eqn\eqAPq{
\ga_{(\ep_1\ldots\ep_m)}(q,t) \eql \prod_{i\in M_- \atop j\in M_+}
  \left(  { \et_{ij}^{(++)}(q,t) \over \et_{ij}^{(-+)}(q,t) } \right) \,,
}
are in the commutant of $S^\pm$.  The coefficients 
$\ga_{(\ep_1\ldots\ep_m)}(q,t)$ satisfy the following properties
\eqn\eqAPna{ 
\ga_{(\ep_1\ldots\ep_m)}(q,1)  \eql 1\,,
}
\eqn\eqAPnb{
\ga_{(\ep_1\ldots\ep_m)}(q^{-1},t^{-1})  \eql 
  \ga_{(\ep_1\ldots\ep_m)}(q,t)\,, 
}
\eqn\eqAPnc{
\ga_{(\ep_1\ldots\ep_m)}(t,q) \eql \ga_{(\ep_1\ldots\ep_m)}(q,t)\,,
}
\eqn\eqAPnd{
\sum_{\ep_i\atop |M_-|=n} \ \ga_{(\ep_1\ldots\ep_m)}(q,t) \eql 
  \sum_{\ep_i\atop |M_-|=m-n} \ \ga_{(\ep_1\ldots\ep_m)}(q,t)\,.
}
Furthermore,
\eqn\eqAPne{
\vth( T_{({\bf a},{\bf b})}(z) ) \eql T_{(-{\bf a},-{\bf b})}(zq^2t^{-2})\,.
}
\medskip

\proof 
Consider the expression \eqAPn, where we normalize $\ga_{(+\ldots+)}=1$.
Using \eqAPc\ it is clear that the $\de(x_iq)$-term in 
$$
[ \Ps_{(\ep_1\ldots\ep_{i-1}+\ep_{i+1}\ldots\ep_m)}(z) ,S^-(w)]
$$
is going to combine with the $\de(x_iqt^{-2})$-term in 
$$
[ \Ps_{(\ep_1\ldots\ep_{i-1}-\ep_{i+1}\ldots\ep_m)}(z) ,S^-(w)]
$$
into a total $t^2$-difference, i.e., that we have
\eqn\eqAPo{
[T_{({\bf a},{\bf b})}(z),S^-(w)] \eql \cD_{t^2}\cdot R(z,w)\,,
}
for some $R(z,w)$, provided we can choose the 
$\ga_{(\ep_1\ldots\ep_m)}(q,t)$ in \eqAPn\ such that
\eqn\eqAPp{
\ga_{(\ep_1\ldots\ep_{i-1}-\ep_{i+1}\ldots\ep_m)}(q,t) \eql
\prod_{j\neq i} \left( { \et_{ij}^{(+\ep_j)}(q,t) 
\over \et_{ij}^{(-\ep_j)}(q,t) }
  \right) \ \ga_{(\ep_1\ldots\ep_{i-1}+\ep_{i+1}\ldots\ep_m)}(q,t) \,,
}
for all choices of $i$ and $\ep_j\,,j\neq i$.
The solution of \eqAPp\ is given by \eqAPq.%
\foot{Note that it is a non-trivial 
fact that the system \eqAPp\ has a solution at all.  Indeed, we have 
$m 2^{m-1}$ equations for $2^m -1$ unknowns.  The dependence of the equations
is however guaranteed by \eqAPm.}
The analysis for the other screening current $S^+(w)$ proceeds similarly.
The properties \eqAPna\ -- \eqAPnc\ trivially follow from 
\eqAPl\ and \eqAPm.  To prove \eqAPnd, consider both sides as a 
(bounded) meromorphic function of one variable, say $z_1$.  Such a function
is uniquely determined by the residues at its poles and its value at 
infinity.  Using the third relation in \eqAPm\ it is easy to check that
the residues of both sides indeed agree.  Finally, \eqAPne\ follows 
straightforwardly from \eqAPnb. \Box\medskip

For parameters $z_i$ in non-generic position, some of the coefficients 
$\ga_{(\ep_1\ldots\ep_m)}(q,t)$ in \eqAPq\ may be vanishing or singular.
Four different situations might occur:
\item{(i)} There exists a pair $(i,j)$ such that $\xi_{ij}=q$.  In this 
case $\et_{ij}^{(++)}(q,t)=0$ such that 
$$
\ga_{(\ep_1\ldots\ep_{i-1}-\ep_{i+1}\ldots\ep_{j-1}+\ep_{j+1}\ldots\ep_m)}
 \eql 0\,,
$$
and the number of terms in \eqAPn\ is effectively reduced by $1/4$.
\item{(ii)} There exists a pair $(i,j)$ such that $\xi_{ij}=t$.  This case
is similar to case (i) due to the symmetry $(q,t)\to(t,q)$.
\item{(iii)} There exists a pair $(i,j)$ such that $\xi_{ij}=1$.  In this 
case $\et_{ij}^{(++)}(q,t)$ and $\et_{ji}^{(++)}(q,t)$, and thus both
$$
\ga_{(\ep_1\ldots\ep_{i-1}-\ep_{i+1}\ldots\ep_{j-1}+\ep_{j+1}\ldots\ep_m)}
$$
and 
$$
\ga_{(\ep_1\ldots\ep_{i-1}+\ep_{i+1}\ldots\ep_{j-1}-\ep_{j+1}\ldots\ep_m)}
$$
are singular.  However, the residue at this singularity vanishes, so the 
expression \eqAPn\ makes perfect sense provided we interpret the right
hand side as the limit of the generic expression in which we let
$\xi_{ij}\to1$.  Note that, by doing so, the $\Ps_{(\ep_1\ldots\ep_m)}(z)$
terms are no longer of the form \eqAPh, but contain derivative terms.
\item{(iv)} There exists a pair $(i,j)$ such that $\xi_{ij}=qt^{-1}$.
In this case $\et_{ij}^{(-+)}(q,t)=0$ and the coefficient 
$$
\ga_{(\ep_1\ldots\ep_{i-1}-\ep_{i+1}\ldots\ep_{j-1}+\ep_{j+1}\ldots\ep_m)}
$$
becomes singular.  This is a genuine singularity.  By renormalizing the
expression for $T_{({\bf a},{\bf b})}(z)$ we obtain an element in the
commutant of $S^\pm$ with terms which are precisely the complement of
the terms remaining in case (i).  Note that, in this case,
$\NO{Y^{(-)}(z_i) Y^{(+)}(z_j)} = 1$.  Thus, the generator obtained 
this way corresponds effectively to an expression \eqAPn\ with $m\to m-2$.
\medskip

Of course, various combinations of the above cases can occur simultaneously.
A particularly important example is when the $z_i$ line up in a single
$q$-string, i.e.,
\eqn\eqAPaa{
z_i \eql zq^{2a+2(i-1)}t^{2b}\,,\qquad i=1,\ldots,m\,,
}
for some $a,b\in\RR$.
One easily verifies that the only nonvanishing coefficients 
$\ga_{(\ep_1\ldots\ep_m)}(q,t)$ in \eqAPn\ are
\eqn\eqAPab{
\ga_{m,n}(q,t) \eqv
\ga_{{{\underbrace{(-\ldots-}\atop n} 
   {\underbrace{+\ldots+)}\atop m-n}}}(q,t)\,,
}
and, therefore, the corresponding element
$T_{({\bf a},{\bf b})}(z)\equiv T_m(z)$ in the commutant
has $m+1$ terms.  Explicitly, $\ga_{m,0}(q,t)=1$ and, for $n=1,\ldots,m$,
\eqn\eqAPac{
\ga_{m,n}(q,t) \eql \prod_{k=1}^n   
{ (q^kt^{-1}-q^{-k}t)(q^{m-k+1} - q^{-(m-k+1)})
 \over (q^k-q^{-k})(q^{m-k+1}t^{-1} - q^{-(m-k+1)}t) } \,.
}
{}From Theorem \thAPa\ it follows that the 
coefficients $\ga_{m,n}(q,t)$ have the following properties
\eqn\eqAPaca{ \eqalign{
\ga_{m,n}(q,1) & \eql 1 \,,\cr
\ga_{m,n}(q^{-1},t^{-1}) & \eql \ga_{m,n}(q,t) \,,\cr
\ga_{m,n}(t,q) & \eql \ga_{m,n}(q,t) \,,\cr
\ga_{m,m-n}(q,t) & \eql \ga_{m,n}(q,t)\,,\cr}
}
which can easily be verified explicitly from \eqAPac.

The foregoing construction of the operators
$T_{({\bf a},{\bf b})}(z)$ is intimately related to the structure of 
the representation ring of the quantum affine algebra
$U_q(\widehat{\frak{sl}_2})$.  Before formulating a more precise conjecture,
let us discuss the extension of the above results to arbitrary simple $\bfg$.
Suppose we start with the operator 
\eqn\eqAPca{
Y_{({\bf a},{\bf b})}(z) \eqv
  \NO{ \prod_{i=1}^\ell\, \prod_{j_i=1}^{m_i} \ Y_i(z_{j_i}^{(i)}) }\,,
}
where 
\eqn\eqAPcb{
z_{j}^{(i)} \eql z q^{2r_i a_{j}^{(i)}} t^{2b_{j}^{(i)}}\,,
}
with $({\bf a},{\bf b})=((a_1^{(1)},b_1^{(1)}),\ldots,(a_{m_\ell}^{(\ell)},
b_{m_\ell}^{(\ell)}))$, 
and try to complete it to an operator in the commutant of all
$S^\pm_i$, i.e., an element of $\cW_{q,t}[\bfg]$.  
In each $\frak{sl}_2$ direction $i$, we may
apply the results of Theorem \thAPa, with the replacement $q\to
q^{r_i}$, and the final result is obtained by pasting together all the
$\frak{sl}_2$ directions.  Clearly, for this algorithm to work,
certain consistency conditions at the intersections of the various
$\frak{sl}_2$ directions must be satisfied.  To illustrate the general  
procedure we have summarized the construction of the generators
of $\cW_{q,t}[\bfg]$ for
the rank $2$ simple Lie algebras $\bfg$ in Appendix A.  

Note that the above algorithm bears close resemblance to the
construction of the irreducible finite dimensional representation
$L(\La)$ of $\bfg$ of highest weight $\La = \sum_i m_i \om_i$.  In
fact, the construction would be exactly identical if, at all
$\frak{sl}_2$ highest weights in direction $i$, the operators $Y_i(z)$
would line up in $q^{r_i}$-strings.  
This is not the case in general,
though, as the completion of $Y_{({\bf a},{\bf b})}(z)$ to an element
in the commutant of $S^\pm_i$ in general has more terms than the
dimension of $L(\La)$.  In fact, the conjecture (cf.\ [\FRc]) is that
the number of terms and their weights are the same as the dimension
and the weights of an irreducible finite dimensional representation
$V$ of the quantum affine algebra $U_q(\whg)$, which decomposes under
$U_q(\bfg)$ as $V\cong L(\La)
\oplus \ldots$, where the dots stands for `smaller' representations.
An example of this, which is worked out in detail in Appendix A.3, is
the generator corresponding to the $\bf 14$ of $U_q(G_2)$, to which an
additional singlet term has to be added, in accordance with the minimal
affinization of the $\bf 14$.  
This conjecture has been verified in all cases where the generators
$T_V(z)$ are explicitly known.  Obviously, when $V$ is one of the fundamental 
representations it coincides with Assumption \thSTad.

The main point here is that the algorithm should hold for other than
the fundamental representations as well [\FRc,\FRe,\FRaff].  This leads to a
generalization of Assumption \thSTad, which will be formulated after
we recall some basic facts from the representation theory of quantum
affine algebras (see, e.g., [\CPd] and references therein).

Let ${\rm Rep\, }(U_q(\whg))$ denote the ring of finite
dimensional representations of $U_q(\whg)$.  
It is well-known that there is a 1--1 correspondence between 
the irreducibles $V \in {\rm Rep\, }(U_q(\whg))$ and monic polynomials
$P_{i,V}(u)$, $i=1,\ldots,\ell$ [\Dr].  Let 
$\{ u_{j_i}^{(i)} \,|\, j_i=1,\ldots,m_i\}$ be the roots of $P_{i,V}(u)$.  
\thm\thAPb
\proclaim Conjecture \thAPb.  Let $V \in {\rm Rep\, }(U_q(\whg))$ be 
irreducible.  We have a map $V\mapsto T_V(z) \in \cW_{q,t}[\bfg]$
given by%
\foot{{\rm Of course, this map is not unique, e.g., it can be
twisted by an automorphism of $U_q(\whg)$.}}
\eqn\eqAPfa{
T_V(z) \eql \sum_{\la\in P(V)} \,
  \sum_{k_\la=1}^{{\rm mult\,}\la}\ c^{V,(k_\la)}_\la(q,t) \,
  Y^{V,(k_\la)}_\la(z)\,,
}
where $\la$ runs over the weights $P(V)$ of $V$, and is uniquely determined
by the above $\frak{sl}_2$-pasting algorithm from the highest weight 
component
\eqn\eqAPfb{
Y_V(z) \eqv Y^V_{\La}(z) \eql \NO{ \prod_{i=1}^\ell \biggl( \prod_{j_i=1}^{m_i}
  \ Y_i(z u_{j_i}^{(i)}) \biggr) }\,,
}
with highest weight $\La = \sum_i m_i \om_i$. \medskip

\ni {\bf Remark.}  In most cases the $Y^{V,(k_\la)}_\la(z)$, with
$\la = \La - \sum_j \al_{i_j}$, will be of 
the form 
\eqn\eqAPfc{
\NO{ Y_V(z)
  A_{i_1}(zq^{a_1}t^{b_1})^{-1}\ldots A_{i_k}(zq^{a_k}t^{b_k})^{-1} }\,.
}
However, it can happen that in some $\frak{sl}_2$ direction there exists
a pair of arguments $(z,w)$ such that $w/z=1$ (case (iii) above), in
which case there will be derivative terms.  This, for example, happens in the
case of the $\cW_{q,t}[G_2]$ generator with highest weight component
$Y_1(zq^{-1})Y_1(zq)$, i.e., the minimal $U_q(\widehat{G_2})$
affinization of the $U_q(G_2)$ irrep.\ $L(2\om_1) \cong {\bf 27}$ which 
is a ${\bf 27} \oplus {\bf 7}$. \medskip

\ni {\bf Remark.} It is known 
that there is a map $V\mapsto t_V(z) = \sum_{m\in\ZZ} t_V[m]z^{-m}$ 
from ${\rm Rep\, }(U_q(\whg))$ 
to generating series of central elements in $U_q(\whg)_{k=-\dcox}$
[\RSTS,\DE], which at $q=1$ reduces to the character of $V$.
These $q$-characters satisfy the following natural properties
\item{(I)} $t_{V\oplus W}(z) = t_V(z) + t_W(z)$, for all $V,W\in 
{\rm Rep\, }(U_q(\whg))$,
\item{(II)} $t_{V\otimes W}(z) = t_V(z)t_W(z)$, for all $V,W\in 
{\rm Rep\,}(U_q(\whg))$,
\item{(III)} $t_{V(a)}(z) = t_V(za)$, for all $V,W\in 
{\rm Rep\, }(U_q(\whg))$, $a\in\CC^*$, and where $V(a)$ is the image of $V$
under the twist automorphism. 

\ni Moreover, it has been shown, at least for $\bfg=\sln$ [\FRa], 
that the evaluation of the image of $t_V(z)$ under the free field
realization $U_q(\whg) \to \cH_{q,1}(\bfg)$ coincides with the
Bethe Ansatz formula for the eigenvalues of the transfer matrix 
corresponding to the finite dimensional representation 
$V \in {\rm Rep\,}(U_q(\whg))$.  
Conjecture \thAPb, which is a slight 
extension of the conjectures in [\FRc,\FRe], 
is a $t$-deformation 
of this classical result in the sense that 
to each $V \in {\rm Rep\,}(U_q(\whg))$ one can associate a field $T_V(z)$
in $\cW_{q,t}[\bfg]$ such that $T_V(z) \to t_V(z)$ for $t\to1$. \medskip

Comparing our explicit calculations in the case of $\frak{sl}_2$ above
with the representation theory of $U_q(\widehat{\sltw})$ [\CPc] shows
that Conjecture \thAPb\ holds for $\sltw$.  The $\sltw$ results
indicate that the tensor product structure of ${\rm Rep\,
}(U_q(\whg))$ is also reflected in $\cW_{q,t}[\bfg]$ through a
quantization of property (II) above (properties (I) and (III) continue
to hold for the quantization $T_V(z)$).  Specifically, since 
 $\cW_{q,t}[\bfg]$ has been defined as a deformed chiral algebra, by the very
axioms of this algebra [\FRb] the residues at the poles of the
(meromorphically continued) composition of two generators are again
elements of $\cW_{q,t}[\bfg]$.  We expect that the composition
$T_V(z) T_W(w)$, $V,W\in {\rm Rep\, }(U_q(\whg))$, in particular
contains poles corresponding to {\it all\/} subquotients of $V\otimes W$.
More precisely\foot{See, also [\FRaff].}
\thm\thAPc
\proclaim Conjecture \thAPc.  Let $V,W \in {\rm Rep\, }(U_q(\whg))$ and 
$T_V(z)$ and $T_W(z)$ be the corresponding elements of $\cW_{q,t}[\bfg]$.
For each subquotient $U$ of $V\otimes W$ there exists a meromorphic 
function $\et_U^{V\otimes W}(x)$ and a choice of $a,b,a',b'\in\ZZ$ such that 
\eqn\eqAPcc{
T_{U}(z q^{a'}t^{b'}) \eql \lim_{w\to zq^at^b}\ \et_U^{VW}({w\over z}) 
  T_V(z) T_W(w)\,.
}\medskip

For $\frak{sl}_2$ the validity of this conjecture,
which is also implicit in [\FRb], follows from the 
observation that (cf.\ [\BP], (2.37)--(2.40))
\eqn\eqAPcd{ \eqalign{
f( {z_2\over z_1} ) Y^{(+)}(z_1) Y^{(-)}(z_2) & \eql 
  { \et^{(++)}_{21}(q,t) \over \et^{(-+)}_{21}(q,t)}
  \NO{ Y^{(+)}(z_1) Y^{(-)}(z_2) } \,, \cr
f( {z_2\over z_1} ) Y^{(-)}(z_1) Y^{(+)}(z_2) & \eql 
{ \et^{(++)}_{12}(q,t) \over \et^{(-+)}_{12}(q,t)} 
  \NO{ Y^{(-)}(z_1) Y^{(+)}(z_2) }\,, \cr}
}
and comparison to the explicit tensor product structure of 
$U_q(\widehat{\frak{sl}_2})$ [\CPc].  In Appendix A we will see
examples for the simple Lie algebras $\bfg$ of rank $2$.

\ni {\bf Remark.} For $V$ and $W$ in generic position the $U_q(\whg)$
module $V\otimes W$ will be irreducible.  In that case we can simply take
$a=b=0$, $\xi=1$ and 
\eqn\eqAPce{
\et_{V\otimes W}^{VW}(x) \eql f_{VW}(x)\,,
}
where $f_{VW}(x)$ is determined by 
\eqn\eqAPcf{
Y_V(z)Y_W(w) \eql f_{VW}({w\over z})^{-1}\NO{Y_V(z)Y_W(w)}\,.
}
For $V\otimes W$ reducible, and $U\subset V\otimes W$ the subquotient 
with highest weight given by the highest weight of $V\otimes W$, the
choice \eqAPce\ suffices as well.  The other subquotients can be projected 
out by using the singularity structure of $(\et_{ij}^{(-+)})^{-1}$ in 
the case of $qt^{-1}$-strings (see case (iv) in the $\frak{sl}_2$ case).

The structure of the commutant for $\frak{sl}_2$ suggests (see, in
particular, \eqAPh) the 
definition of a `deformed Weyl group action' as follows.  Define, for each
$i=1,\ldots,\ell$, a transformation $T_i: \cH_{q,t}[\bfg] 
\to \cH_{q,t}[\bfg]$ by
\eqn\eqAPba{
T_i\, y_j[m] \eql y_j[m] - \de_{ij}\, p_i^{m}\, a_i[m]\,,
}
or, equivalently, 
\eqn\eqAPbb{
T_i\, a_j[m] \eql a_j[m] - C_{ij}(q^m,t^m)\, p_i^{m}\, 
a_i[m]\,,
}
where $p_i \equiv q^{r_i}t^{-1}$.
Let us denote by $T_{q,t}[\bfg]$ the algebra generated by the $T_i$,
$i=1,\ldots,\ell$.  It is easy to check that $T_{q,t}[\bfg]$ acts by
automorphisms of $\cH_{q,t}[\bfg]$.  The $T_i$'s satisfy various
relations, the simplest of which are
\eqn\eqAPbc{
(T_i - 1) (T_i + p_i^{-2d}) \eql 0 \,,
}
or, equivalently,
\eqn\eqAPbd{
T_i^{-1} \eql p_i^{2d} T_i + (1- p_i^{2d}) \,,
}
where $d$ is defined in \STTxb.  In addition,
\eqn\eqAPbe{ \eqalign{
T_i T_j & \eql  T_j T_i\,, \qquad {\rm if} \quad    C_{ij}C_{ji}\eql0\,, \cr
T_i T_j T_i & \eql  T_j T_iT_j \,, \qquad {\rm if}\quad  
  C_{ij}C_{ji}\eql1\,.\cr}
}
Note that, for $\frak{sl}_2$, we have 
\eqn\eqAPbf{
\Ps_{(-\ldots-)}(z) \eql T \Ps_{(+\ldots+)}(z) \,.
}
Thus, the algebra $T_{q,t}[\bfg]$ can be used to construct `Weyl orbits' of 
terms in the expressions for the commutant.

In particular, let $w_0=r_{i_1}\ldots r_{i_n}$, ($n=|\De_+|$),
be a reduced expression
for the longest Weyl group element $w_0$.  Define 
\eqn\eqAPbg{
T_{w_0} \eql T_{i_1} \ldots T_{i_n}\,,
}
then, as one can verify on a case by case basis,%
\foot{Albeit not completely obvious, the result is independent of the 
choice of reduced expression for $w_0$.}
\eqn\eqAPbh{
T_{w_0} Y_i(z) \eql Y_{i^*}(z q^{-r^\vee\dcox}t^{\cox})^{-1}\,.
}

Finally, let us return to the proofs of Theorems \thSTTfx\ and \thSTae\
in Section 2.  Under Assumption \thSTad, Theorem \thSTTfx\ 
follows immediately from \eqAPna\ -- \eqAPnd, while Theorem \thSTae\
follows from \eqAPbh\ by using, in particular, the assumption
that $T_i(z)$ is uniquely determined by its highest component.
Moreover, both theorems generalize to $\cW_{q,t}[\bfg]$ elements with
more general highest weight component \eqAPca.  The generalization of
Theorem \thSTTfx\ is self-evident, while the generalization of Theorem
\thSTae\ is (cf.\ \eqAPne\  for $\bfg = \sltw$)
\eqn\eqAPbi{
\vth( T_{ ({\bf a},{\bf b})} (z) ) \eql 
  T_{ (-{\bf a}^*,-{\bf b}^*)} (zq^{r^\vee\dcox}t^{-\cox})\,,
}
where
\eqn\eqAPbj{
(a_j^{*(i)},b^{*(i)}_j) \eqv (a^{(i^*)}_j,b^{(i^*)}_j)\,.
}


\newsec{Involutions, contravariant forms and the Kac determinant}

The algebra $\cH_{q,t}'[\bfg]$ has a natural two-parameter family of
anti-involutions $\om_{a,b}$, $a,b\in\CC$, defined by
\eqn\STTaa{\eqalign{
\om_{a,b}(a_i[0])&\eql -a_i[0]+2r^\vee(\rh,\al_i)-2\be(\rh^\vee,\al_i)\,,\cr
\om_{a,b}(a_i[m])&\eql -(q^at^b)^{m} a_i[-m]\,,\qquad m\not=0\,,\cr}}
and
\eqn\STTax{
\om_{a,b}(Q_i)\eql Q_i\,,}
or, equivalently,
\eqn\STTab{
\om_{a,b}(A_i(z))\eql A_i({q^at^b\over z})^{-1}\,,\qquad 
\om_{a,b}(Y_i(z))\eql Y_i({q^at^b\over z})^{-1}\,.}
Also
\eqn\STTac{
\om_{a,b}(S_i^\pm(z))\eql {1\over z^2}S_i^{\pm}({q^at^b\over z})\,.
} 
In particular, the screening operators $S_i^\pm$ are invariant
under $\om_{a,b}$, up to a multiplicative factor, and therefore
$\om_{a,b}$ is well-defined on $\cW_{q,t}[\bfg]$.
\thm\thSTTa
\proclaim Lemma \thSTTa. The action of $\om_{a,b}$ on the generators
$T_i(z)$ of $\cW_{q,t}[\bfg]$ is given by
\eqn\STTad{
\om_{a,b}(T_i(z))\eql T_{i^*}({q^{a+r^\vee \dcox}t^{b-\cox}\over z})\,.}
\medskip

\proof Note that the action of $\om_{a,b}$ on $T_i(z)$ in \STWbe\ is
obtained by composing the transformation $\vth$ of \STWbf\ with the
transformation $z\rightarrow q^at^b/z$.  The lemma follows by using
Theorem \thSTae\ and then applying the transformation $z\rightarrow
q^at^b/z$.  \Box\medskip

In the following we will set $a=-r^\vee \dcox$ and $b=\cox$, and
denote the corresponding anti-involution by $\om$. Obviously, we then
have $\om(T_i(z))=T_{i^*}(1/z)$.

The anti-involution $\om$ induces a unique contravariant bilinear form
$\langle- |-\rangle_F$ on $F(\bmu)\times F(\mu)$ such that
$\langle\bmu|\mu\rangle_F=1$, where $\bmu=-\mu+2r^\vee\rh-2\be\rh^\vee$
(cf.\ \TSSac).
We will denote by $g_{\sbfla\sbfla'}=\langle\bmu |\om(y[-\bfla])
y[-\bfla']|\mu\rangle_F$ the matrix elements of this
form in the basis \STTxc.  Similarly, $\om$ induces a unique 
contravariant bilinear form $\langle- |-\rangle$
on $M(h)\times M(h)$ with matrix elements 
$G_{\sbfla\sbfla'}=\langle h |\om(T[-\bfla])T[-\bfla']|h\rangle$.
Clearly, both $g_{\sbfla\sbfla'}$ and $G_{\sbfla\sbfla'}$
vanish unless $|\bfla|=|\bfla'|$.  The map $\imath:M(h(\mu)) \to
F(\mu)$ is an isometry for generic values of $(q,t;q^\mu)$.
\thm\thSTTba
\proclaim Lemma \thSTTba. Let $g^{(n)}(q,t)=
\det( g_{\sbfla\sbfla'} )_{|\sbfla|=|\sbfla'|=n}$ be the determinant 
at level $n$. Then
\eqn\STTbe{g^{(n)}(q,t) \eql C_n \prod_{r,s\geq1\atop rs\leq n} \Biggl(
(q^{r^\vee \dcox} t^{-\cox})^{{\ell r}} \biggl( {\rm det}\, M(q^r,t^r)\biggr) 
(q^r-q^{-r} )^\ell (t^r-t^{-r} )^\ell 
\Biggr)^{p_\ell(n-rs)} \,,
}
where $p_\ell(m)$ is the number of multi-partitions of $m$ 
and $C_n$ is a constant independent of $q$ and $t$.
\medskip

We give an explicit formula for ${\rm det}\,C(q,t)$ and 
${\rm det}\,D(q,t)$ for all simple
Lie algebras $\bfg$, from which 
${\rm det}\,M(q,t)$ follows. \medskip

\settabs 8\columns
\+ $A_\ell$. & ${\rm det}\ C = ({p^{\ell+1 } - p^{-\ell-1})/(p - p^{-1}})$ \cr
\+ $B_\ell$. & ${\rm det}\ C = q^{2\ell-1}t^{-\ell} + q^{-2\ell+1}t^{\ell}$, 
   &&&&${\rm det}\ D = (q+q^{-1})^{\ell-1}$ \cr
\+ $C_\ell$. & ${\rm det}\ C = q^{\ell+1}t^{-\ell} + q^{-\ell-1}t^{\ell}$,
   &&&& ${\rm det}\ D = q+q^{-1}$ \cr
\+ $D_\ell$. & ${\rm det}\ C = (p+p^{-1}) ( p^{\ell-1} + p^{-\ell+1})$ \cr
\+ $E_6$. & ${\rm det}\ C  = p^6 + p^4 -1 + p^{-4} + p^{-6}$ \cr
\+ $E_7$. & ${\rm det}\ C  = p^7 + p^5 - p - p^{-1} + p^{-5} + p^{-7}$\cr
\+ $E_8$. & ${\rm det}\ C = p^8 + p^6 - p^2 - p^{-2} + p^{-6} + p^{-8}$ \cr
\+ $F_4$. &  ${\rm det}\ C = q^6t^{-4} -1 + q^{-6}t^4$, &&&&
   ${\rm det}\ D = (q+q^{-1})^2$ \cr
\+ $G_2$. & ${\rm det}\ C = q^4t^{-2} -1 + q^{-4}t^2$, &&&&
   ${\rm det}\ D = q^2 +1 + q^{-2}$ \cr
\medskip

For $\bfg$ simply-laced the expressions above follow from the lemma below,
the expressions for $B_\ell$ and $C_\ell$ were given in [\FRc].  The 
remaining cases were computed by brute force.
\thm\thTSab
\proclaim Lemma \thTSab.  The eigenvalues $\la_i(q,t)$ of $C_{ij}(q,t)$ 
for $\bfg$ simply-laced are given by
\eqn\TSSba{
\la_i(q,t) \eql (p+p^{-1}) + 2 \cos
\left( {\pi e_i\over \cox} \right)\,,
}
where $p=qt^{-1}$, and $e_i$, $i=1,\ldots,\ell$, are the exponents of $\bfg$.
\medskip

Let $\al\in\De$, $\mu\in\bfh^*$, and $r,s\in\NN$.  Define
\eqn\TSSbb{
G_\al^{(r,s)}(q,t;q^\mu) \eql 
q^{-r^\vee (\rh,\al) +  {1\over2} r^\vee s (\al,\al) } t^{(\rh^\vee,\al) -r }
q^{(\mu,\al)} - 
q^{r^\vee (\rh,\al) - {1\over2} r^\vee s (\al,\al) } t^{-(\rh^\vee,\al) +r }
q^{-(\mu,\al)}\,,
}
then, under the action of $W$ (see \TSSaa), we have 
\eqn\eqTSSbc{
G_\al^{(r,s)}(q,t;q^{w*\mu}) \eql G_{w^{-1}\al}^{(r,s)}(q,t;q^\mu)
}
In addition,
\eqn\TSSbd{
G_\al^{(r,s)}(q,t;q^{\bmu}) \eql
G_{-\al}^{(r,s)}(q,t;q^\mu)\,,
}
where 
\eqn\TSSbda{
\bmu \eql -\mu + 2r^\vee \rh -2 \be \rh^\vee\,.
}

We define the matrix of $\imath$ in the standard bases \STTxh\ and
\STTxc\  by
\eqn\STTxxa{
\imath \left(T[-\bfla]|h(\mu)\rangle\right)
\eql\sum_{\sbfla'}\Pi_{\sbfla\sbfla'}(q,t;q^\mu)\, y[-\bfla']|\mu\rangle\,,
}
and denote by $\Pi^{(n)}(q,t;q^\mu)\equiv \det
(\Pi_{\sbfla\sbfla'})_{|\sbfla|=|\sbfla'|=n}$ the determinant of this
matrix at level $n$.
\thm\thTSSba
\proclaim Theorem \thTSSba.  Given Assumption \thSTad,
the Kac determinant of $\cW_{q,t}[\bfg]$ at level $n$
is given by
\eqn\TSSbe{
G^{(n)}(q,t;q^\mu) \eql 
\Pi^{(n)}(q,t;q^{\bmu})\
g^{(n)}(q,t) \ \Pi^{(n)}(q,t;q^{\mu})\,,
}
where 
\eqn\TSSbf{
\Pi^{(n)}(q,t;q^{\mu}) \eql C_n \prod_{r,s\geq1\atop rs\leq n} 
\Biggl( (q^{r^\vee \dcox} t^{-\cox})^{-{\ell r\over2}}\ 
\prod_{\al\in\De_+} G_\al^{(r,s)}(q,t;q^\mu) \Biggr)^{p_\ell(n-rs)}\,,
}
and $g^{(n)}(q,t)$ is given in \STTbe.
That is, using \TSSbd,
\eqn\TSSbg{
G^{(n)}(q,t;q^{\mu})  \eql C_n \prod_{r,s\geq1\atop rs\leq n} 
\Biggl( \biggl( \prod_{\al\in\De} G_\al^{(r,s)}(q,t;q^{\mu})  \biggr)
\biggl({\rm det}\, M(q^r,t^r)\biggr) (q^r-q^{-r} )^\ell (t^r-t^{-r} )^\ell 
\Biggr)^{p_\ell(n-rs)}\,.
}
\medskip

\proof The proof is a generalization of the proof in [\BP] for $\bfg =
\frak{sl}_2$ to which we refer for more details.  First we observe that 
the determinant $G^{(n)}(q,t;q^{\mu})$ can be factorized as in \TSSbe\ 
by using
the norm preserving homomorphism $\imath:M(h(\mu)) \to F(\mu)$.  The
Fock space determinant $g^{(n)}(q,t)$ was computed in Lemma \thSTTba,
so it remains to compute a sufficient number of vanishing lines of
$\Pi^{(n)}(q,t;q^\mu)$ (note that $\Pi^{(n)}(q,t;q^\mu)$ is a Laurent 
polynomial in $q^{\mu}$).
The construction of a set of vanishing lines of $\Pi^{(n)}(q,t;q^\mu)$ 
proceeds as follows.  For every weight $\mu$ of the form 
$\mu=\be \mu^{(+)} - r^\vee\mu^{(-)}$ with 
$\mu^{(+)}\in P_+^\vee$ and $\mu^{(-)}\in P_+$ and $i=1,\ldots,\ell$, we
can construct a $\cW_{q,t}[\bfg]$ singular vector in $F(\bmu)$ 
(where $\bmu$ is given by \TSSbda) at level $(\mu^{(+)}+\rh^\vee,\al_i)
(\mu^{(-)}+\rh,\al_i^\vee)$.  Explicitly, this singular vector is the 
image of the highest weight vector $|\bmu + r \be \al_i^\vee\rangle$ under
the map 
\eqn\TSSya{
\oint \prod_{j=1}^r \left( dz_jS_+^i(z_j) \right) \ : \ 
F(\bmu + r \be \al_i^\vee) ~\rightarrow~ F(\bmu)\,,
}
where 
\eqn\TSSyb{ 
r \eql (\mu^{(+)}+\rh^\vee,\al_i) \,, \qquad
s \eql (\mu^{(-)}+\rh,\al_i^\vee)\,.
}
Note that, with the definition \TSSyb,
\eqn\TSSyc{
(\mu - r^\vee \rh + \be \rh^\vee, \al_i) \eql r \be - \half r^\vee s (\al_i,
  \al_i) \,.
}
For $\cW_{q,t}[\sln]$ this construction was carried out in [\AKOSa],
where it was also shown that in this case the result could be
expressed in terms of Macdonald polynomials.  The construction for
general $\bfg$ is a straightforward generalization (cf.\ [\BS] in the
conformal case).  Due to the non-degenerate pairing between $F(\mu)$
and $F(\bmu)$ (for generic values of $(q,t;q^\mu)$) there must exist a
vector in the cokernel of the map $\imath : M(h(\mu)) \to F(\mu)$ at
level $rs$, i.e., by using \TSSyc, we conclude that
$\Pi^{(n)}(q,t;q^\mu)$ has vanishing lines
$G_{\al_i}^{(r,s)}(q,t;q^\mu)$ for all $i=1,\ldots,\ell$, and $rs\geq n$.
Using the Weyl group invariance \eqTSSbc, one then proves that 
$\Pi^{(n)}(q,t;q^\mu)$ is given by \TSSbf\ up to 
a Laurent polynomial $C_n$ in $q$, $t$ and $q^\mu$.  To prove that 
$C_n$ is actually a constant it suffices to compute the leading order 
term in $\Pi^{(n)}(q,t;q^{\mu})$ (the (partial) ordering is given by $q^\mu 
\succeq q^{\mu'}$ iff $\mu-\mu' \in \ZZ_{\geq 0}\cdot \De_+$).
We find (cf.\ [\BP] for more details)
\eqn\TSSbh{\eqalign{
\prod_{\sbfla \vdash n} &  \ \prod_{i=1}^\ell \biggl(
q^{-2r^\vee(\rh,\om_i)} t^{2(\rh^\vee,\om_i)} q^{2(\mu,\om_i)} \biggr)
^{{\rm length}(\la^{(i)})} \cr & \eql
\prod_{r,s\geq1\atop rs\leq n}\ \prod_{\al\in\De_+} \ \biggl( 
q^{-r^\vee(\rh,\al)} t^{(\rh^\vee,\al)} q^{(\mu,\al)} \biggr)^{p_\ell(n-rs)}\cr
& \eql \prod_{r,s\geq1\atop rs\leq n}\ \Biggl( 
(q^{r^\vee \dcox} t^{-\cox})^{-{\ell r\over2}} \prod_{\al\in\De_+} \ \biggl(
q^{-r^\vee (\rh,\al) +  {1\over2}r^\vee s (\al,\al) } t^{(\rh^\vee,\al) -r }
q^{(\mu,\al)}\biggr) \Biggr)^{p_\ell(n-rs)} \,,\cr}
}
where we have used
\eqn\TSSbi{
\sum_i \om_i \eql \half \sum_{\al\in\De_+} \al\,,
}
and 
\eqn\TSSbj{
\ell\, \cox \eql 2|\De_+| \,,\qquad 
 \ell \, \dcox \eql \sum_{\al\in\De_+} (\al,\al)\,.
}
This concludes the proof of the theorem. \Box\medskip

\noindent {\bf Remark.} It follows from Theorem \thSTae\ that the 
determinant of the matrix
\eqn\TSSbk{
\Pi^{\prime(n)}_{\sbfla\sbfla'} (q,t;q^\mu) \eql 
(q^{r^\vee \dcox} t^{-\cox})^{{n\over2}}\ 
\Pi^{(n)}_{\sbfla\sbfla'} (q,t;q^\mu)\,,
}
satisfies (cf.\ (3.19) in [\BP])
\eqn\TSSbl{
\Pi^{\prime(n)}(q^{-1},t^{-1};q^{-\mu}) \eql
\Pi^{\prime(n)}(q,t;q^\mu) \,.
}
This can indeed be verified from 
the explicit expression \TSSbf\ of this determinant. \medskip


\newsec{The center of $\cW_{q,t}[\bfg]$ at roots of unity}

In this section we consider the limit of $\cW_{q,t}[\bfg]$ when $t^2$
is  a primitive $\nroot$-th root of unity, $t^2\rightarrow
\sqrtN1$, and $q^2$ is generic, i.e., we have $t^{{2\nroot}}=1$ and
$t^{2j}\not=1$ for all $j=1,\dots,{\nroot}-1$, while $q^{2j}\not=1$
for all $j\not=0$. We will restrict our discussion mainly to the
simplest case $\bfg=\slN$. Note that, because of the duality
$\cW_{q,t}[\bfg]\simeq \cW_{t,q}[\bfg]$, which holds when $\bfg$ is
simply laced, it does not matter with respect to which deformation
parameter the limit is taken. For $\bfg$ non-simply laced the
situation is different as the deformation is not symmetric under the
interchange of $q$ and $t$. (In fact, the duality above is then
replaced by a more complicated relation [\FRc].) It follows from the
definitions in Section 2 that the dependence of $\cW_{q,t}[\bfg]$ on
$t$, unlike that on $q$, is in some sense universal for all $\bfg$.
This allows an extension of the construction to the general case,
which is then verified on examples for rank 2 algebras using the
explicit realizations from Appendix A.

The definition of $\cW_{q,t}[\bfg]$, in the limit $t^2\rightarrow
\sqrtN{1}$, presents a subtlety in that the oscillators \STWam\ and
thus the screening currents \STWal\ are not well defined in this limit. 
Note, however, that the divergent factor in
the definition of the oscillators \STWam\ cancels out in the
commutators between the screening operators $S_i^-$ and the fields in
${\bf H}_{q,t}[\bfg]$. Thus the problem of computing $\cW_{q,t}[\bfg]$
as a commutant defined through those equations is well posed, also in the
limit, and will simply lead to an algebra $\cW_{q,t}[\bfg]$ with generators
$T_i(z)$, as in Assumption \thSTad, specialized from the generic case to 
$t^2=\sqrtN1$.

It is an obvious observation that, as follows from
\STWaa\ and \STWag,  in the limit $t^2\rightarrow \sqrtN1$ (or 
$q^2\rightarrow \sqrtN1$) the algebra $\cH_{q,t}[\bfg]$ has a large
center generated by the oscillators $a_i[m]$ (equivalently, $y_i[m]$),
with $m=0\,{\rm mod}\,{\nroot}$. In turn, this implies that there
should be a corresponding center of $\cW_{q,t}[\bfg]$ and our goal is
to verify that by constructing this center explicitly in terms of the
generators $T_i(z)$.

The existence of this center may also be inferred from the formula for
the Kac determinant of $\cW_{q,t}[\bfg]$ in \TSSbg. Namely, the
determinant, $G^{(n)}$, contains a factor
$\prod_r(q^r-q^{-r})(t^r-t^{-r})$, and thus vanishes for either $t^2$
or $q^2$ a root of unity and $n$ sufficiently large.  Thus, for those
values of the deformation parameters, the Verma module should have
additional singular vectors that are independent of the highest
weight. Indeed, it follows from the nonvanishing of \TSSbf\ for a
generic $q$ and $\mu$, that the Verma module, $M(h(\mu))$, is
isomorphic with the Fock module, $F(\mu)$. Obviously, the latter has
infinitely many singular vectors corresponding to the center of
$\cH_{q,t}[\bfg]$, which in turn give rise to singular vectors in the
Verma module independent of the highest weight.

To make this discussion more explicit let us now  consider the case
$\bfg=\sln$. The generators of $\cW_{q,t}[\bfg]$ are given
in Theorem \thSTac. Note that if we recast \STWbb\  as in
\STWbe, all the coefficients $c^{\om_i}_\la(q,t)$ are equal to one and 
\eqn\STFab{Y^{\om_i}_\la(z)\eql 
\NO{\La_{l_1}(zp^{-i+1})\La_{l_2}(zp^{-i+3})\ldots \La_{l_i}(zp^{i-1})}\,. 
}
The (nondegenerate) weight $\la\in P(V(\om_i))$ in \STFab\
corresponds to the sequence $(l_1,\ldots,l_i)$,
$l_1<\ldots<l_i$, such that  $\la=\sum_j \ep_{l_j}$ where
$\ep_i\,,i=1,\ldots,N$, is an overcomplete basis in terms of which the
simple roots of $\sln$ are given as $\al_i=\ep_i-\ep_{i+1}$.  In
particular, $(1,2,\ldots,i)$ is the sequence for $\om_i$, in
accordance with
\eqn\STFbxx{
Y^{\om_i}_{\om_i}(z)\eql Y_i(z) \eql
\NO{ \La_1(zp^{-i+1})\La_{2}(zp^{-i+3})\ldots \La_i(zp^{i-1}) }\,.
}
{}For $\la\in P(V(\om_i))$ and $\la'\in P(V(\om_j))$, define
$f_{\la\la'}^{ij}(x)$ by
\eqn\STFba{
Y_{\la}^{\om_i}(z)Y_{\la'}^{\om_j}(w)\eql f_{\la\la'}^{ij}
({w\over z})^{-1}  \NO{Y_{\la}^{\om_i}(z)Y_{\la'}^{\om_j}(w) }\,.
}
Setting $\la=\la'=\om_i$, we obtain $f_{\om_i\om_i}^{ii}(x)=
f_{ii}(x)$, which is given in \STWda. Then an arbitrary
$f_{\la\la'}^{ij}(x)$ can be computed using [\FF]:
\eqn\STFlac{
\La_l(z)\La_{l'}(w) \eql s_{l,l'}({w\over z}) f_{11}({w\over z})^{-1}
\NO{ \La_l(z)\La_{l'}(w)}\,,
}
where
\eqn\STFlab{
s_{l,l'}(x)\eql
\cases{	s(x)&for $\quad l<l'\,,$\cr
1&for $\quad l=l'\,,$\cr
	s(xp^2)&for $\quad l>l'\,,$\cr}\qquad 
}
and
\eqn\STFssa{
s(x)\eql {(q-q^{-1}x)(t^{-1}-tx)\over (p-p^{-1}x)(1-x)}\,.
}

For our purposes it will suffice to consider only weights in the same
representation, i.e., we take $\la,\la'\in P(V(\om_i))$ corresponding
to sequences $\la=(l_1,\ldots,l_i)$ and $\la'=(l'_1,\ldots,l'_i)$,
respectively.

An immediate consequence of \STFlab\ is that $f_{\la\la}^{ii}(x)$ does
not depend on a particular choice of the weight,
\eqn\STFbb{
f_{\la\la}^{ii}(x)\eql f_{\om_i\om_i}^{ii}(x)\eql f_{ii}(x)\,.
}
{}For $\la\not=\la'$, we can use \STFlab\ to compute the additional
factor that arises from the points at which $(l_1,\ldots,l_i)$ and
$(l'_1,\ldots,l'_i)$ are different. Hence we write
\eqn\STFcc{
f_{\la\la'}^{ii}(x)^{-1}\eql 
s_{\la\la'}(x) f_{ii}(x)^{-1}\,,
}
where
\eqn\STFdd{
s_{\la\la'}(x)~\equiv~s_{l_1\ldots l_i,l'_1\ldots l'_i}(x)\eql 
\prod_{a=1}^i s_{l_a,l'_1\ldots l'_i}(x)\,,
}
and $s_{l_a,l'_1\ldots l'_i}(x)$ is the contribution due to  the
ordering of $l_a$ with respect to $l'_1,\ldots,l'_i$. We find (cf.\
[\FF])
\eqn\STFee{
s_{l_a,l'_1\ldots l'_i}(x)\eql\cases{
1& for  $\quad l_a=l'_b\,,$\cr
s(xp^{2m})& for  $\quad l'_{a+m-1}<l_a<l'_{a+m}\,,$\cr}
}
where the (in)equality is satisfied for some $1\leq b\leq m$ or
$-a\leq m\leq i-a+1$.

The center can  be constructed in terms of the generators $T_i(z)$
by generalizing the corresponding result for $\sltw$ in [\BP].

\thm\thSFaa
\proclaim Theorem \thSFaa. For $t^2=\sqrtN1$ and $q$ generic, define 
\eqn\STFxvc{
\Psi_i^{(k)}(z)\eql 
\lim_{z_m\rightarrow zt^{2(\nroot-m)} }\biggl(\prod_{m<n}
f_{ii}({z_n\over z_m})
\biggr) T_{i}(z_1)\ldots T_i(z_\nroot)\,.}
Then we have
\eqn\STFac{
\Psi_i^{(k)}(z) \eql
\sum_{\la\in P(V(\om_i))}\NO{ Y^{\om_i}_\la(zt^{2(\nroot-1)})
Y^{\om_i}_\la(zt^{2(\nroot-2)})\ldots 
Y^{\om_i}_\la(z) }\,.
}
and $\Psi_i^{(k)}(z)$ is a well defined series of central elements of
$\cW_{q,t}[\sln]$.
\medskip

\proof After expanding  \STFxvc\ using \STWbb, \STFab\
and \STFba\ we obtain a sum of terms of the form
\eqn\STFax{
\biggl(\prod_{m<n} s_{\la_m\la_n}({z_n\over z_m})\biggr)
\NO{ Y_{\la_1}^{\om_i}(z_1)\ldots Y_{\la_\nroot}^{\om_i}(z_\nroot)}\,.
}
It follows from \STFssa\ that for a generic $q$, and thus $p$, none of
the factors $s_{\la_m\la_n}(z_n/z_m)$ develop a pole in the limit 
$z_n/z_m \rightarrow t^{2(m-n)}$. (Note that we have $|m-n|<\nroot$.) 

Consider the weights $\la_m=(l_{m,1},\ldots,l_{m,i})$ and
$\la_{m+1}=(l_{m+1,1},\ldots,l_{m+1,i})$.  Using \STFee\ it is easy to see
that $s_{\la_m\la_{m+1}}(t^{-2})$ has a vanishing factor of $s(t^{-2})$
unless $l_{m,a}\geq l_{m+1,a}$ for all $a=1,\ldots,i$.  Next consider
$\la_1$ and $\la_\nroot$. By the previous argument we may assume
$\la_{1,a}\geq \la_{2,a}\geq\ldots\geq \la_{\nroot,a}$. Suppose
$\la_{1,i}>\la_{\nroot,i}$. This results in a factor 
$s(p^2t^{2\nroot-2})$, which vanishes for $t^2=\sqrtN1$. Thus we must
have $l_{1,i}=l_{2,i}=\ldots=l_{\nroot,i}$. Proceeding by induction we
then find that the only nonvanishing terms in \STFax\ arise for
$\la_1=\ldots=\la_\nroot$, which proves \STFac.

It follows from the explicit expressions \STWah, \STWai\ and \STWbc,
that for $t^2=\sqrtN1$ all terms in the sum on the right hand side of \STFac\
have an expansion in terms of the oscillators $a_i[n\nroot]$,
$n\in\ZZ$. Thus, at least formally, \STFac\ is in the center of
$\cW_{q,t}[\bfg]$. \Box
\medskip

If one tries to expand the right hand side of \STFxvc\ in terms of the
modes $T_i[n]$, the resulting series is divergent. It is however well
defined when acting on the vacuum of a Verma module. We then obtain
a series of singular vectors that are manifestly independent of the
highest weight. We refer the reader to [\BP] for explicit formulae for
the low lying singular vectors in the $\sltw$ case. 

Here let us consider as an example the case of $\slth$ with $k=2$,
i.e., $t^2=-1$. For a generic $q$ and $h$ we find the following pairs
of singular vectors at levels $2m$, $m\geq 1$:
\eqn\STFyxa{\eqalign{
\Ps_1^{(2)}[-2m] |h\rangle  
\eql \biggl( \sum_{\sbfla\vdash 2m} m_{\sbfla}(1,t^2)
f_{11}&(R_{12})  T_1[-\la_1] T_1[-\la_2]\cr& - 2 (-1)^m p^{2m} 
[4m-3]_p\,
T_2[-2m] \biggr) 
|h\rangle \,,\cr
\Ps_2^{(2)}[-2m] |h\rangle \eql 
\biggl( \sum_{\sbfla\vdash 2m} m_{\sbfla}(1,t^2)
f_{11}&(R_{12}) T_2[-\la_1] T_2[-\la_2] \cr&- 2 (-1)^m 
p^{-2m} 
[4m-3]_p\,
T_1[-2m] \biggr) |h\rangle \,,\cr} } 
where the sum runs over partitions $\bfla=(\la_1,\la_2)$,
$\la_1\geq\la_2\geq0$ of $2m$, $R_{12}$ is the raising operator acting by
\eqn\STFyxz{
R_{12}T_i[m]T_i[n]\eql T_i[m-1]T_i[n+1]\,, } 
and $m_{\sbfla}$ is the monomial symmetric function [\MD]. We can make these
formulae even more explicit for $m=1$ where we find
\eqn\STFzxc{\eqalign{
\Ps_1^{(2)}[-2] |h\rangle & \eql -\biggl( T_1[-1]T_1[-1]  + 
{2\over [3]_p} 
T_1[-2]T_1[0] - {2 p^2} \,T_2[-2] \biggr) |h\rangle\,, \cr
\Ps_2^{(2)}[-2] |h\rangle & \eql -\biggl( T_2[-1]T_2[-1]  +  
{2\over [3]_p} 
T_2[-2]T_2[0] - {2\over p^2} \,T_1[-2] \biggr) |h\rangle \,.\cr}
}

The above discussion has a simple generalization to arbitrary
$\cW_{q,t}[\bfg]$.

\thm\centconjj
\proclaim Conjecture \centconjj. For $t^2=\sqrtN1$ and $q$ generic,
define $\Psi_i^{(k)}(z)$ as in \STFxvc. Then
\eqn\STFacx{
\Psi_i^{(k)}(z) \eql
\sum_{\la\in P(V(\om_i))}\sum_{j_\la=1}^{{\rm mult}\,\la}
 \NO{ Y^{\om_i,(j_\la)}_\la(zt^{2\nroot-2})
Y^{\om_i,(j_\la)}_\la(zt^{2\nroot-4})\ldots 
Y^{\om_i,(j_\la)}_\la(z) }\,,
}
is a series of central elements in $\cW_{q,t}[\bfg]$.
\medskip

The last assertion in the conjecture is obvious, provided we prove the
expansion \STFacx. We have verified Conjecture \centconjj\ in all cases where
the generators $T_i(z)$ are known explicitly. 

As remarked before, for simply laced $\bfg$, the center for $q^2$  
a primitive $k$-th root of unity (and $t$ generic) follows from 
Conjecture \centconjj\ by using the duality invariance $(q,t) \to (t,q)$.
For non-simply laced, the situation $q^2= \sqrtN{1}$ is more complicated.
In particular, the generating series of central elements will in general
not be homogeneous of fixed order in the generators $T_i(z)$, due to a
different rescaling $q\to q^{r_i}$ in the various $\sltw$ directions.  We
will leave this case for further investigation.

\bigskip
\noindent {\bf Acknowledgements.} P.B.\ is supported by a \qeii\ research 
fellowship from the Australian Research Council and K.P.\ is supported
in part by the U.S.\ Department of Energy Contract
\#DE-FG03-84ER-40168.
\bigskip

\vfil\eject

\appendix{A}{Examples -- the $\cW_{q,t}[\bfg]$ algebras of rank $2$}

In this appendix we illustrate some of the ideas of the paper in the case
of the deformed $\cW$-algebras $\cW_{q,t}[\bfg]$ corresponding to the 
rank $2$ simple Lie algebras $A_2$, $B_2$ and $G_2$.  We provide explicit
expressions for the generators and their relations and illuminate the 
connection to the representation theory of the quantum affine algebra
$U_q(\whg)$.

To simplify the notation, let us define
\eqn\eqGBb{
\cst{a_1,\ldots,a_r}{b_1,\ldots,b_s} \eql
  { (a_1-a_1^{-1}) \ldots (a_r-a_r^{-1}) \over (b_1 - b_1^{-1})\ldots
 (b_s - b_s^{-1}) }\,.
}
Also, we recall that the function $f^{ij}_{\la\la'}(x)$ is defined 
by (cf.\ \STFba)
\eqn\eqGBa{
Y^{\om_i}_\la(z) Y^{\om_j}_{\la'}(w) \eql
  f^{ij}_{\la\la'}( {w\over z})^{-1} 
  \NO{ Y^{\om_i}_\la(z) Y^{\om_j}_{\la'}(w) } \,.
}

\appsubsec{$\cW_{q,t}[\frak{sl}_3]$}

The case $\bfg = A_2 = \frak{sl}_3$ has been discussed in detail in
[\FF,\AKOSa].  For completeness we give a brief review.

We adopt the following conventions.  In terms of an
overcomplete basis of $\RR^2$ given by vectors $\{ \ep_1,
\ep_2,\ep_3\}$, satisfying 
\eqn\eqGGba{
\ep_i \cdot \ep_j \eql \de_{ij} - \textstyle{1\over3}\,,\qquad 
\ep_1+\ep_2+\ep_3 \eql 0\,,
}
the simple roots and fundamental weights of $\frak{sl}_3$ are written 
as 
\eqn\eqGGbb{ \eqalign{
\al_1 \eql \ep_1 - \ep_2 \,, & \qquad \al_2 \eql \ep_2 - \ep_3 \,,\cr
\om_1 \eql \ep_1 \,, & \qquad \om_2 \eql \ep_1 + \ep_2 \eql -\ep_3\,.\cr}
}
We have $(r_1,r_2)=(1,1)$, and $r^\vee=1$, $\cox=\dcox =3$.
The weights of the irreducible $\frak{sl}_3$ representations 
$L(\om_1) = {\bf 3}$ and $L(\om_2)= {\bf 3}^*$ are given by
$\{\ep_1,\ep_2,\ep_3\}$ and $\{-\ep_1,-\ep_2,-\ep_3\}$, respectively. 

The deformed Cartan matrix is given by
\eqn\eqGGbc{
C_{ij}(q,t) \eql \left( \matrix{  qt^{-1} + q^{-1}t & -1 \cr
                                  -1 & qt^{-1} + q^{-1}t \cr} \right)\,.
}
In particular,
\eqn\eqGGbd{ \eqalign{
A_1(z) & \eql \NO{ Y_1(zqt^{-1}) Y_1(zq^{-1}t) Y_2(z)^{-1} } \,,\cr
A_1(z) & \eql \NO{ Y_1(z)^{-1} Y_2(zqt^{-1}) Y_2(zq^{-1}t) } \,.\cr}
}

The generators of $\cW_{q,t}[\frak{sl}_3]$ follow by applying Theorem \thAPa\
in the various $\sltw$ directions.  One finds  (see also \STWbc)
\eqn\eqGGad{ 
T_i (z) \eql \sum_{\la\in P(V(\om_i))} \ c_\la^{\om_i}(q,t)\, 
  Y^{\om_i}_\la(z) \,,
}
where
\eqn\eqGGaa{ \eqalign{
\La_1(z) & \eqv Y^{\om_1}_{\ep_1}(z)  \eql \NO{ Y_1(z) } \,,\cr
\La_2(z) & \eqv Y^{\om_1}_{\ep_2}(z)
           \eql \NO{ Y^{\om_1}_{\ep_1}(z)A_1(zq^{-1}t)^{-1} } 
           \eql \NO{ Y_1(zq^{-2}t^2)^{-1} Y_2(zq^{-1}t) } \,,\cr
\La_3(z) & \eqv Y^{\om_1}_{\ep_3}(z)
           \eql \NO{ Y^{\om_1}_{\ep_2}(z) A_2(zq^{-2}t^2)^{-1} } 
           \eql \NO{ Y_2(zq^{-3}t^3)^{-1} } \,,\cr}
}
\eqn\eqGGab{ \eqalign{
Y^{\om_2}_{-\ep_3}(z) & \eql \NO{ Y_2(z) } \,,\cr
Y^{\om_2}_{-\ep_2}(z) & \eql \NO{ Y^{\om_2}_{-\ep_3}(z)A_2(zq^{-1}t)^{-1} }
                        \eql \NO{ Y_1(zq^{-1}t)Y_2(zq^{-2}t^2)^{-1} } \,,\cr
Y^{\om_2}_{-\ep_1}(z) & \eql \NO{ Y^{\om_2}_{-\ep_2}(z)A_1(zq^{-2}t^2)^{-1} } 
                        \eql \NO{ Y_1(zq^{-3}t^3)^{-1} } \,.\cr}
}
and all $c_\la^{\om_i}(q,t)=1$, in agreement with \STTxk.  

Note that we can write (cf.\ \STWbb)
\eqn\eqGGac{ \eqalign{
Y^{\om_2}_{-\ep_3}(z) & \eql \NO{ \La_1(zq^{-1}t)\La_2(zqt^{-1}) } \,,\cr
Y^{\om_2}_{-\ep_2}(z) & \eql \NO{ \La_1(zq^{-1}t)\La_3(zqt^{-1}) } \,,\cr
Y^{\om_2}_{-\ep_1}(z) & \eql \NO{ \La_2(zq^{-1}t)\La_3(zqt^{-1}) } \,.\cr}
}
In fact, it is not hard to see that, 
\eqn\eqGGae{
T_2 (z qt^{-1}) \eql \lim_{w\to zq^2t^{-2}} \ 
  f_{\ep_1\ep_2}^{11}( {w\over z} ) \, T_1(z) T_1(w)\,,
}
which illustrates Conjecture \thAPc.  Similarly, one can verify, e.g.,
\eqn\eqGGaf{ \eqalign{
T_{V(2\om_1)}(zq) & \eql \lim_{w\to zq^2} \ 
  f_{\ep_1\ep_1}^{11}( {w\over z} ) \, T_1(z) T_1(w)\,,\cr
T_{V(\om_1+\om_2)}(z) & \eql T_1(z)T_2(zq) \,,\cr}
}
where $T_{V(2\om_1)}(z)$ and $T_{V(\om_1+\om_2)}(z)$ are the 
$\cW_{q,t}[\frak{sl}_3]$ generators corresponding to the
$V(2\om_1) = {\bf 6}$ and $V(\om_1+\om_2)={\bf 8}$ of 
$U_q(\widehat{\frak{sl}_3})$, respectively.

For the commutation relations one finds
\eqn\eqGGa{ \eqalign{
f_{11}({w\over z}) T_1(z) T_1(w)  - & f_{11}({z\over w}) T_1(w) T_1(z) \cr 
  \eql & \cst{q,t^{-1}}{qt^{-1}} 
  \biggl( \de( q^2t^{-2}{w\over z}) T_2(wqt^{-1}) -
  \de( q^{-2}t^{2}{w\over z}) T_2(zqt^{-1}) \biggr)\,, \cr
 & \cr
f_{12}({w\over z}) T_1(z) T_2(w)  - & f_{12}({z\over w}) T_2(w) T_1(z) \cr 
  \eql & \cst{q,t^{-1}}{qt^{-1}} \biggl( \de( q^3t^{-3}{w\over z}) 
  - \de( q^{-3}t^{3}{w\over z}) \biggr) \,,\cr
 & \cr
f_{22}({w\over z}) T_2(z) T_2(w)  - & f_{22}({z\over w}) T_2(w) T_2(z) \cr 
  \eql &  \cst{q,t^{-1}}{qt^{-1}} 
  \biggl( \de( q^2t^{-2}{w\over z}) T_1(wqt^{-1}) -
  \de( q^{-2}t^{2}{w\over z}) T_1(zqt^{-1}) \biggr)\,, \cr}
}
where $f_{ij}(x)$ is defined in \STFbx.
Note that $f_{11}(x)=f_{22}(x)$.

\appsubsec{$\cW_{q,t}[B_2]$}

In this appendix we compute explicit expressions 
for the generators and relations of the deformed $\cW_{q,t}[B_2]$ algebra. 
The classical limit $t\to1$, i.e., $\cW_{q,1}[B_2]$,
has been discussed already in [\FRa,\FRc].

We adopt the following conventions.  In terms of an orthonormal basis 
$\{\ep_1,\ep_2\}$ of $\RR^2$, the simple roots and 
fundamental weights of $B_2$ are given by
\eqn\eqGBab{ \eqalign{
\al_1 \eql \ep_1 - \ep_2 \,, & \qquad \al_2 \eql \ep_2 \,,\cr
\om_1 \eql \ep_1 \,,& \qquad \om_2 \eql \half (\ep_1 +\ep_2) \,.\cr}
}
We have $(r_1,r_2)=(2,1)$ and $r^\vee=2$, $\cox=4$, $\dcox=3$.
The weights of the $B_2$ irreducible representations $L(\om_1) = {\bf 5}$
and $L(\om_2)= {\bf 4}$ are given by $\{ \pm\ep_1, \pm\ep_2, 0\}$ and
$\{ \half (\pm\ep_1 \pm\ep_2) \}$, respectively.

The deformed Cartan matix is given by
\eqn\eqGBaa{
C_{ij}(q,t) \eql \left( \matrix{ q^2t^{-1} + q^{-2}t & -1 \cr
                         -(q+q^{-1}) & qt^{-1} + q^{-1}t \cr} \right)\,.
}
In particular,
\eqn\eqGBzz{ \eqalign{
A_1(z) & \eql \NO{ Y_1(zq^2t^{-1}) Y_1(zq^{-2}t) 
  Y_2(zq^{-1})^{-1}Y_2(zq)^{-1} } \,,\cr
A_2(z) & \eql \NO{ Y_1(z)^{-1} Y_2(zqt^{-1}) Y_1(zq^{-1}t) } \,.\cr}
}

Using the $\sltw$ pasting procedure outlined in Section 3 we find the
following generators of $\cW_{q,t}[B_2]$:%
\foot{In fact, $T_2(z)$ of $\cW_{q,t}[B_2]$ coincides with $T_1(z)$ of
$\cW_{q,t}[C_2]$, so that those generators can be read-off from the
results in [\FRc].}
\eqn\eqBbh{
T_i(z) \eql \sum_{\la \in P(V(\om_i))} \ c^{\om_i}_\la(q,t)\ 
  Y^{\om_i}_\la(z)\,,
}
where
\eqn\eqBbb{ \eqalign{
\La_1(z) & \eqv Y^{\om_1}_{\ep_1}(z)  
           \eql \NO{ Y_1(z) } \,,\cr
\La_2(z) & \eqv Y^{\om_1}_{\ep_2}(z) 
           \eql \NO{ Y^{\om_1}_{\ep_1}(z) A_1(zq^{-2}t)^{-1}} 
           \eql \NO{ Y_1(zq^{-4}t^2)^{-1} Y_2(zq^{-3}t)Y_2(zq^{-1}t)}\,,\cr
\La_0(z) & \eqv Y^{\om_1}_{0}(z) 
           \eql \NO{ Y^{\om_1}_{\ep_2}(z) A_2(zq^{-4}t^2)^{-1} }
           \eql \NO{ Y_2(zq^{-5}t^3)^{-1} Y_2(zq^{-1}t) } \,,\cr
\La_{\bar{2}}(z) & \eqv Y^{\om_1}_{-\ep_2}(z) 
           \eql \NO{ Y^{\om_1}_{0}(z) A_2(zq^{-2}t^2)^{-1} }
           \eql \NO{ Y_1(zq^{-2}t^2)Y_2(zq^{-5}t^3)^{-1}Y_2(zq^{-3}t^3)^{-1}} 
           \,,\cr
\La_{\bar{1}}(z) & \eqv Y^{\om_1}_{-\ep_1}(z) 
           \eql \NO{ Y^{\om_1}_{-\ep_2}(z)A_1(zq^{-4}t^4)^{-1} }
           \eql  \NO{ Y_1(zq^{-6}t^4)^{-1} }\,,\cr}
}
\eqn\eqBba{ \eqalign{
Y^{\om_2}_{ {1\over2} (\ep_1+\ep_2)}(z) & \eql \NO{ Y_2(z) } \,,\cr
Y^{\om_2}_{ {1\over2} (\ep_1-\ep_2)}(z) & 
     \eql \NO{ Y^{\om_2}_{ {1\over2} (\ep_1+\ep_2)}(z)A_2(zq^{-1}t)^{-1} }
     \eql \NO{ Y_1(zq^{-1}t) Y_2(zq^{-2}t^2)^{-1} } \,, \cr
Y^{\om_2}_{ - {1\over2} (\ep_1-\ep_2)}(z) & 
     \eql \NO{ Y^{\om_2}_{ {1\over2} (\ep_1-\ep_2)}(z) A_1(zq^{-3}t^2)^{-1}}
     \eql \NO{ Y_1(zq^{-5}t^3)^{-1}  Y_2(zq^{-4}t^2)} \,, \cr
Y^{\om_2}_{ - {1\over2} (\ep_1+\ep_2)}(z) & 
     \eql \NO{ Y^{\om_2}_{ - {1\over2} (\ep_1-\ep_2)}(z) A_2(zq^{-5}t^3)^{-1}}
     \eql \NO{ Y_2(zq^{-6}t^4)^{-1} } \,,\cr}
}
and the coefficients $c^{\om_i}_\la(q,t)$ are given by
\eqn\eqBbi{
c^{\om_1}_0(q,t) \eql \ga_{2,1}(q,t) \eql \cst{q^2,qt^{-1}}{q,q^2t^{-1}}
  \,,\qquad c^{\om_1}_\la(q,t) \eql 1\,,\quad \la\neq0\,.
}
\eqn\eqBbg{
c^{\om_2}_\la(q,t) \eql 1\,,\qquad \forall \la\,.
}
The construction of $T_1(z)$ illustrates an important feature that 
does not occur in the $\cW_{q,t}[\sln]$ case, namely the occurance of
an $\sltw$ string with 3 terms (the one built on $Y^{\om_1}_{\ep_2}(z)$),
and consequently the nontrivial coefficient $c^{\om_1}_0(q,t) 
=\ga_{2,1}(q,t)$.  The generators above are in perfect agreement
with Assumption \thSTad\ as well as Theorems \thSTTfx\ and \thSTae.

Note that we can write
\eqn\eqBbc{ \eqalign{
\La_1(z) & \eql \NO{ Y^{\om_2}_{ {1\over2} (\ep_1+\ep_2)}(zq^{-1}t) 
  Y^{\om_2}_{ {1\over2} (\ep_1-\ep_2)}(zqt^{-1}) }\,,\cr 
\La_2(z) & \eql \NO{ Y^{\om_2}_{ {1\over2} (\ep_1+\ep_2)}(zq^{-1}t) 
  Y^{\om_2}_{ - {1\over2} (\ep_1-\ep_2)}(zqt^{-1})} \,,\cr 
\La_0(z) & \eql \NO{ Y^{\om_2}_{ {1\over2} (\ep_1+\ep_2)}(zq^{-1}t) 
  Y^{\om_2}_{ - {1\over2} (\ep_1+\ep_2)}(zqt^{-1})} \,,\cr 
\La_{\bar{2}}(z) & \eql \NO{ Y^{\om_2}_{ {1\over2} (\ep_1-\ep_2)}(zq^{-1}t) 
  Y^{\om_2}_{ - {1\over2} (\ep_1+\ep_2)}(zqt^{-1})} \,,\cr 
\La_{\bar{1}}(z) & \eql \NO{ Y^{\om_2}_{ - {1\over2} (\ep_1-\ep_2)}(zq^{-1}t) 
  Y^{\om_2}_{ - {1\over2} (\ep_1+\ep_2)}(zqt^{-1})} \,.\cr}
}
In fact, we have
\eqn\eqBbca{
T_1(z qt^{-1}) \eql \lim_{w\to zq^2t^{-2} } \ f_{ {1\over2}(\ep_1+\ep_2),
  {1\over2}(\ep_1-\ep_2)}^{11}( {w\over z} ) T_2(z)T_2(w)\,,
}
once more illustrating Conjecture \thAPc.

For the commutation relations one finds 
\eqn\eqGBc{ \eqalign{
f_{11}({w\over z}) T_1(z) T_1(w)  - & f_{11}({z\over w}) T_1(w) T_1(z) \cr 
  \eql & \cst{q^2,t^{-1}}{q^2t^{-1}} 
  \biggl( \de( q^4t^{-2}{w\over z}) T_{V(2\om_2)}(wq^2t^{-1}) -
  \de( q^{-4}t^{2}{w\over z}) T_{V(2\om_2)}(zq^2t^{-1}) \biggr) \cr
 & + \cst{q^2,t^{-1},q^3t^{-1},qt^{-2}}{qt^{-1},q^2t^{-1},q^3t^{-2}}
 \biggl( \de( q^6t^{-4}{w\over z}) -
  \de( q^{-6}t^{4}{w\over z}) \biggr) \,,\cr
 & \cr
f_{12}({w\over z}) T_1(z) T_2(w)  - & f_{12}({z\over w}) T_2(w) T_1(z) \cr 
  \eql & \cst{q^2,t^{-1}}{q^2t^{-1}} \biggl( \de( q^5t^{-3}{w\over z}) 
  T_2(zq^{-1}t) - \de( q^{-5}t^{3}{w\over z}) T_2(zqt^{-1})\biggr) \,,\cr
 & \cr
f_{22}({w\over z}) T_2(z) T_2(w)  - & f_{22}({z\over w}) T_2(w) T_2(z) \cr 
  \eql &
  \cst{q,t^{-1}}{qt^{-1}} \biggl( \de( q^2t^{-2}{w\over z}) T_1(wqt^{-1}) -
  \de( q^{-2}t^{2}{w\over z}) T_1(zqt^{-1}) \biggr) \cr
 & + \cst{q,t^{-1},q^3t^{-1},q^2t^{-2}}{qt^{-1},q^2t^{-1},q^3t^{-2}}
 \biggl( \de( q^6t^{-4}{w\over z}) -
  \de( q^{-6}t^{4}{w\over z}) \biggr) \,,
\cr}
}
where
\eqn\eqGBd{ \eqalign{
T_{V(2\om_2)}(z) & \eql \lim_{w\to zq^2}
 f_{22}({w\over z}) T_2(zq^{-1})T_2(wq^{-1}) \cr 
  & \eql \NO{ Y_2(zq^{-1})Y_2(zq) } + \ldots\,,\cr}
}
is the $\cW_{q,t}[B_2]$ generator corresponding to the 
irreducible $U_q(\widehat{B_2})$ representation that decomposes under 
$U_q(B_2)$ as a ${\bf 10} \oplus {\bf 1}$.


\appsubsec{$\cW_{q,t}[G_2]$}

In this appendix we illustrate our algorithm to compute explicit expressions 
for the generators and relations of the algebra $\cW_{q,t}[G_2]$.  
For $t=1$, i.e.\ $\cW_{q,1}[G_2]$, the generators and (part of) the 
Poisson algebra structure were already discussed in [\Kog].

We adopt the following conventions.  The simple roots are normalized as 
\eqn\eqAPda{
(\al_1,\al_1) \eql \textstyle{2\over3}\,,\qquad (\al_2,\al_2) \eql 2\,,
\qquad (\al_1,\al_2) \eql -1\,,
}
such that $r_1=1$, $r_2=3$ and $r^\vee=3$.  The Coxeter and dual Coxeter
numbers are $\cox=6$, $\dcox=4$, while the deformed Cartan matrix is
given by
\eqn\eqAPdb{
C_{ij}(q,t) \eql \left( \matrix{  qt^{-1}+q^{-1}t & -(q^2+1+q^{-2}) \cr 
   -1 & q^3t^{-1} + q^{-3}t \cr} \right) \,,
}
such that, in particular,
\eqn\eqAPdc {\eqalign{
A_1(z) & \eql \NO{ Y_1(zqt^{-1}) Y_1(zq^{-1}t) Y_2(z)^{-1} } \,,\cr
A_2(z) & \eql \NO{ Y_1(zq^2)^{-1}Y_1(z)^{-1}Y_1(zq^{-2})^{-1}Y_2(zq^3t^{-1})
 Y_2(zq^{-3}t) }\,. \cr}
}

The two fundamental representations of $G_2$ are $L(\om_1) = {\bf 7}$
and $L(\om_2) = {\bf 14}$.  While the representation $L(\om_1)$ can
be affinized to a 
finite dimensional $U_q(\widehat{G_2})$ module $V(\om_1)$, the minimal 
affinization $V(\om_2)$ of $L(\om_2)$ involves the addition of a 
singlet [\Dr], i.e., as a representation
of $U_q(G_2)$ this $V(\om_2)$ decomposes as ${\bf 14} \oplus {\bf 1}$.

The corresponding $\cW_{q,t}[G_2]$ generators are given by
\eqn\eqAPdd{
T_i(z) \eql \sum_{\la\in P(V(\om_i))} \ c^{\om_i}_\la(q,t)\, 
  Y^{\om_i}_\la(z)\,,
}
where
\eqn\eqAPca{ \eqalign{
\La_1(z) & \eqv Y^{\om_1}_{2\al_1+\al_2}(z) \eql \NO{ Y_1(z) } \,,\cr
\La_2(z) & \eqv Y^{\om_1}_{\al_1+\al_2}(z) 
           \eql \NO{ \La_1(z) A_1(z q^{-1}t)^{-1} }
           \eql \NO{ Y_1(zq^{-2} t^2)^{-1} Y_2(zq^{-1} t) }\,,\cr
\La_3(z) & \eqv Y^{\om_1}_{\al_1}(z) 
           \eql \NO{ \La_2(z) A_2(zq^{-4}t^{2})^{-1} }
           \eql \NO{ Y_1(zq^{-6} t^2) Y_1(zq^{-4} t^2)
                     Y_2(zq^{-7} t^3 )^{-1}}\,,\cr
\La_4(z) & \eqv Y^{\om_1}_{0}(z) 
           \eql \NO{ \La_3(z) A_1(zq^{-7}t^{3})^{-1}}
           \eql \NO{ Y_1(zq^{-8} t^4)^{-1} Y_1(zq^{-4} t^2) } \,,\cr
\La_5(z) & \eqv Y^{\om_1}_{-\al_1}(z) 
           \eql \NO{ \La_4(z) A_1(zq^{-5}t^{3})^{-1}}
           \eql \NO{ Y_1(zq^{-8} t^4)^{-1} Y_1(zq^{-6} t^4)^{-1}
                     Y_2(zq^{-5} t^3)}\,,\cr
\La_6(z) & \eqv Y^{\om_1}_{-(\al_1+\al_2)}(z) 
           \eql \NO{ \La_5(z) A_2(zq^{-8}t^{4})^{-1} }
           \eql \NO{ Y_1(zq^{-10} t^{4})Y_2(zq^{-11} t^5)^{-1}  }\,,\cr
\La_7(z) & \eqv Y^{\om_1}_{-(2\al_1+\al_2)}(z) 
           \eql \NO{ \La_6(z) A_1(zq^{-11}t^{5})^{-1}}
           \eql \NO{Y_1(zq^{-12} t^{6})^{-1} }\,,\cr}
}
\vfil\eject
while 
\eqn\eqAPeb{ \eqalign{
Y^{\om_2}_{3\al_1+ 2\al_2}(z) 
        & \eql \NO{Y_2(z)} \,,\cr
Y^{\om_2}_{3\al_1 + \al_2}(z)  
        & \eql \NO{Y^{\om_2}_{3\al_1+ 2\al_2}(z)A_2(zq^{-3}t)^{-1}} \cr
        & \eql \NO{ Y_1(zq^{-5}t)Y_1(zq^{-3}t)Y_1(zq^{-1}t)
               Y_2(zq^{-6}t^2)^{-1}}\,,\cr
Y^{\om_2}_{2\al_1 + \al_2}(z)
        & \eql\NO{Y^{\om_2}_{3\al_1 + \al_2}(z) A_1(zq^{-6}t^2)^{-1}}\cr
	& \eql \NO{ Y_1(zq^{-3}t)Y_1(zq^{-1}t)Y_1(zq^{-7}t^3)^{-1}}\,,\cr
Y^{\om_2}_{\al_1 + \al_2}(z)
        & \eql\NO{Y^{\om_2}_{2\al_1 + \al_2}(z)A_1(zq^{-4}t^2)^{-1}}\cr
	& \eql \NO{ Y_1(zq^{-1}t)Y_1(zq^{-7}t^3)^{-1}Y_1(zq^{-5}t^3)^{-1}
		Y_2(zq^{-4}t^2)} \,,\cr
Y^{\om_2}_{\al_1}(z)
        & \eql\NO{ Y^{\om_2}_{\al_1 + \al_2}(z) A_2(zq^{-7}t^3)^{-1}} \cr
        & \eql\NO{Y_1(zq^{-9}t^3) Y_1(zq^{-1}t) Y_2(zq^{-10}t^4)^{-1} }\,,\cr
Y^{\om_2}_{\al_2}(z) 
        & \eql \NO{Y^{\om_2}_{\al_1 + \al_2}(z)A_1(zq^{-2}t^2)^{-1}}\cr
        & \eql \NO{ Y_1(zq^{-7}t^3)^{-1} Y_1(zq^{-5}t^3)^{-1} 
          Y_1(zq^{-3}t^3)^{-1}Y_2(zq^{-4}t^2) Y_2(zq^{-2}t^2) } \,,\cr 
Y^{\om_2}_0(z) & \eql \NO{ Y^{\om_2}_{\al_1}(z)A_1(zq^{-10}t^4)^{-1}} \cr
        & \eql \NO{Y_1(zq^{-11}t^5)^{-1} Y_1(zq^{-1}t) }\,,\cr
Y^{\om_2\prime}_0(z) & \eql \NO{ Y^{\om_2}_{\al_1}(z)A_1(zq^{-2}t^2)^{-1}} 
          \eql \NO{ Y^{\om_2}_{\al_2}(z) A_2(zq^{-7}t^3)^{-1} } \cr
        & \eql \NO{ Y_1(zq^{-9}t^3)Y_1(zq^{-3}t^3)^{-1}Y_2(zq^{-10}t^4)^{-1}
          Y_2(zq^{-2}t^2) }\,,\cr
Y^{\om_2\prime\prime}_0(z)  
        & \eql \NO{ Y^{\om_2}_{\al_2}(z) A_2(zq^{-5}t^3)^{-1} } \cr
        & \eql \NO{ Y_2(zq^{-8}t^4)^{-1}Y_2(zq^{-4}t^2) }\,, \cr
Y^{\om_2}_{-\al_2}(z)
	& \eql\NO{Y^{\om_2 \prime}_{0}(z)A_2(zq^{-5}t^3)^{-1}  }\cr
	& \eql \NO{Y_1(zq^{-9}t^3)Y_1(zq^{-7}t^3)Y_1(zq^{-5}t^3)
		Y_2(zq^{-10}t^4)^{-1}Y_2(zq^{-8}t^4)^{-1} }\,,\cr
Y^{\om_2}_{-\al_1}(z)
	& \eql\NO{Y^{\om_2 \prime}_0(z) A_1(zq^{-10}t^4)^{-1} }\cr
	& \eql\NO{Y_1(zq^{-11}t^5)^{-1}Y_1(zq^{-3}y^3)^{-1}
			Y_2(zq^{-2}t^2) }\,,\cr
Y^{\om_2}_{-\al_1-\al_2}(z)
	& \eql\NO{Y^{\om_2}_{-\al_1}(z)A_2(zq^{-5}t^3)^{-1} }
	  \eql\NO{Y^{\om_2}_{-\al_2}(z)A_1(zq^{-10}t^4)^{-1} }\cr
	& \eql\NO{Y_1(zq^{-11}t^5)^{-1}Y_1(zq^{-7}t^3)
		Y_1(zq^{-5}t^3)Y_2(zq^{-8}t^4)^{-1} }\,,\cr
Y^{\om_2}_{-2\al_1-\al_2}(z)
	& \eql\NO{Y^{\om_2}_{-\al_1-\al_2}(z)A_1(zq^{-8}t^4)^{-1} }\cr
	& \eql\NO{Y_1(zq^{-11}t^5)^{-1}Y_1(zq^{-9}t^5)^{-1}
	Y_1(zq^{-5}t^3)}\,,\cr
}}
\vfil\eject
$$ \eqalign{
Y^{\om_2}_{-3\al_1-\al_2}(z)
	& \eql\NO{Y^{\om_2}_{-2\al_1-\al_2}(z)A_1(zq^{-6}t^4)^{-1} }\cr
	& \eql\NO{Y_1(zq^{-11}t^5)^{-1}Y_1(zq^{-9}t^5)^{-1}
	Y_1(zq^{-7}t^5)^{-1}Y_2(zq^{-6}t^4)}\,,\cr
Y^{\om_2}_{-3\al_1-2\al_2}(z)
	& \eql\NO{Y^{\om_2}_{-3\al_1-2\al_2}(z)A_2(zq^{-9}t^5)^{-1} }\cr
	& \eql\NO{Y_2(zq^{-12}t^6)^{-1}} \,,\cr
} $$
and the coefficients $c^{\om_i}_\la(q,t)$ are given by 
\eqn\eqAPde{ \eqalign{
c^{\om_1}_{2\al_1+\al_2}(q,t) & \eql c^{\om_1}_{\al_1+\al_2}(q,t) 
  \eql c^{\om_1}_{\al_1}(q,t) \eql c^{\om_1}_{-\al_1}(q,t) \cr
& \eql c^{\om_1}_{-(\al_1+\al_2)}(q,t) 
  \eql c^{\om_1}_{-(2\al_1+\al_2)}(q,t) \eql 1\,,\cr
c^{\om_1}_0(q,t) & 
  \eql \ga_{2,1}(q,t) \eql \cst{q^2,qt^{-1}}{q,q^2t^{-1}}\,. \cr}
}
and 
\eqn\eqAPec{ \eqalign{
c^{\om_2}_{3\al_1+2\al_2}(q,t) & \eql c^{\om_2}_{3\al_1+\al_2}(q,t)
  \eql c^{\om_2}_{\al_2}(q,t)\eql c^{\om_2}_{-\al_2}(q,t)\cr 
& \eql c^{\om_2}_{-3\al_1-\al_2}(q,t)
  \eql c^{\om_2}_{-3\al_1-2\al_2}(q,t)\eql 1\,,\cr
  c^{\om_2}_{2\al_1+\al_2}(q,t) & \eql c^{\om_2}_{\al_1+\al_2}(q,t) \eql 
  c^{\om_2}_{\al_1}(q,t) \eql c^{\om_2}_{-\al_1}(q,t) \cr
& \eql c^{\om_2}_{-\al_1-\al_2}(q,t) 
  \eql c^{\om_2}_{-2\al_1-\al_2}(q,t) \eql 
  \cst{q,q^3t^{-1}}{q,q^3t^{-1}} \,,\cr
c^{\om_2}_0(q,t) & \eql \cst{q^3,qt^{-1},q^5t^{-1},q^4t^{-2}}
  {q,q^3t^{-1},q^4t^{-1},q^5t^{-2}}\,,\cr
c_0^{\om_2\prime}(q,t) & \eql \cst{q^4,qt^{-1}}{q,q^4t^{-1}} \,,\cr
c_0^{\om_2\prime\prime}(q,t) &\eql -\cst{q^2,qt}{q,q^2t^{-1}} \,.\cr}
}\medskip

Note that again we find perfect agreement with both Theorem \thSTTfx\
and \thSTae.  The construction of the generator $T_2(z)$ illustrates
two important features.  First, the $\sltw$ string built on, e.g.,
$Y^{\om_2}_{\al_1}(z)$, requires $4$ terms as compared to the
$3$-dimensional $U_q(\sltw)$ representation which occurs at this point
in the ${\bf 14}$ of $U_q(G_2)$.  This illustrates the necessity for
extending the ${\bf 14}$ by a ${\bf 1}$.  Secondly, the $\sltw$
strings built on $Y^{\om_2}_{\al_1}(z)$ and $Y^{\om_2}_{\al_2}(z)$
intersect at the point $Y^{\om_2\prime }_0(z)$.  For consistency of
the $\sltw$ pasting procedure we therefore need to find the same
coefficient $c^{\om_2\prime}_0(q,t)$ regardless of which $\sltw$ path
we choose to reach $Y^{\om_2\prime }_0(z)$.  This can indeed be
verified.

Note that 
\eqn\eqAPega{ \eqalign{
Y^{\om_2}_{3\al_1+ 2\al_2}(z) & \eql \NO{ \La_1(zq^{-1}t)\La_2(zqt^{-1})}\,,\cr
Y^{\om_2}_{3\al_1 + \al_2}(z) & \eql \NO{ \La_1(zq^{-1}t)\La_3(zqt^{-1})}\,,\cr
Y^{\om_2}_{2\al_1 + \al_2}(z) & \eql \NO{ \La_1(zq^{-1}t)\La_4(zqt^{-1})}\,,\cr
Y^{\om_2}_{\al_1 + \al_2}(z)  & \eql \NO{ \La_1(zq^{-1}t)\La_5(zqt^{-1})}\,,\cr
Y^{\om_2}_{\al_1}(z)          & \eql \NO{ \La_1(zq^{-1}t)\La_6(zqt^{-1})}\,,\cr
Y^{\om_2}_{\al_2}(z)          & \eql \NO{ \La_2(zq^{-1}t)\La_5(zqt^{-1})}\,,\cr
Y^{\om_2}_0(z)                & \eql \NO{ \La_1(zq^{-1}t)\La_7(zqt^{-1})}\,,\cr
Y^{\om_2\prime}_0(z)          & \eql \NO{ \La_2(zq^{-1}t)\La_6(zqt^{-1})}\,,\cr
Y^{\om_2\prime\prime}_0(z)    & \eql \NO{ \La_3(zq^{-1}t)\La_5(zqt^{-1})}\,,\cr
Y^{\om_2}_{-\al_2}(z)         & \eql \NO{ \La_3(zq^{-1}t)\La_6(zqt^{-1})}\,,\cr
Y^{\om_2}_{-\al_1}(z)         & \eql \NO{ \La_2(zq^{-1}t)\La_7(zqt^{-1})}\,,\cr
Y^{\om_2}_{-\al_1-\al_2}(z)   & \eql \NO{ \La_3(zq^{-1}t)\La_7(zqt^{-1})}\,,\cr
Y^{\om_2}_{-2\al_1-\al_2}(z)  & \eql \NO{ \La_4(zq^{-1}t)\La_7(zqt^{-1})}\,,\cr
Y^{\om_2}_{-3\al_1-\al_2}(z)  & \eql \NO{ \La_5(zq^{-1}t)\La_7(zqt^{-1})}\,,\cr
Y^{\om_2}_{-3\al_1-2\al_2}(z) & \eql \NO{ \La_6(zq^{-1}t)\La_7(zqt^{-1})}\,,\cr
}}
so that we can also write
\eqn\eqAPeg{ 
T_2(z) \eql  \sum_{i<j} \ c_{ij}(q,t) \NO{\La_i(zq^{-1}t)\La_j(zqt^{-1})}\,,
}
with appropriately chosen coefficients $c_{ij}(q,t)$ (some of which are
vanishing).  In fact, an explicit examination of all the contractions
shows
\eqn\eqAPei{
T_2 (z qt^{-1}) \eql \lim_{w\to zq^2t^{-2}} \ 
  f_{2\al_1+\al_2,\al_1+\al_2}^{11}( {w\over z} ) \, T_1(z) T_1(w)\,,
}
which again confirms Conjecture \thAPc.
\vfil\eject

The commutation relations are given by
\eqn\eqBGa{ \eqalign{
f_{11}({w\over z})&  T_1(z) T_1(w)  -  f_{11}({z\over w}) T_1(w) T_1(z) \cr 
  \eql &
  \cst{q,t^{-1}}{qt^{-1}} \biggl( \de( q^2t^{-2}{w\over z}) T_2(wqt^{-1}) -
  \de( q^{-2}t^{2}{w\over z}) T_2(zqt^{-1}) \biggr) \cr
 & + \cst{q^2,t^{-1},q^4t^{-1},q^3t^{-2}}{q^2t^{-1},q^3t^{-1},q^4t^{-2}}
  \biggl( \de( q^8t^{-4}{w\over z}) T_1(wq^4t^{-2}) -
  \de( q^{-8}t^{4}{w\over z}) T_1(zq^4t^{-2}) \biggr) \cr
 & + \cst{q,t^{-1},q^4t^{-1},q^2t^{-2},q^6t^{-2},q^5t^{-3}}
  {qt^{-1},q^2t^{-1},q^4t^{-2},q^5t^{-2},q^6t^{-3}}
 \biggl( \de( q^{12}t^{-6}{w\over z}) -
  \de( q^{-12}t^{6}{w\over z}) \biggr) \,,\cr
 & \cr
f_{12}({w\over z})& T_1(z) T_2(w)  -  f_{12}({z\over w}) T_2(w) T_1(z) \cr 
  \eql & \cst{q^3,t^{-1}}{q^3t^{-1}}
  \biggl( \de(q^7t^{-3} {w\over z}) T_{V(2\om_1)}(wq^2t^{-1}) -
  \de(  q^{-7}t^{3}{w\over z}) T_{V(2\om_1)}(wq^{-2}t) \biggr) \cr
 & + \cst{q^3,t^{-1},q^5t^{-1},q^2t^{-2}}{q^2t^{-1},q^3t^{-1},q^5t^{-2}}
  \biggl( \de(q^{11}t^{-5} {w\over z}) T_1(wqt^{-1}) -
  \de( q^{-11}t^{5}{w\over z}) T_1(wq^{-1}t) \biggr)\,, \cr
 & \cr
f_{22}({w\over z})& T_2(z) T_2(w)  -  f_{22}({z\over w}) T_2(w) T_2(z) \cr 
 \eql & \cst{q^3,t^{-1}}{q^3t^{-1}} 
  \biggl( \de( q^6t^{-2}{w\over z}) T_{V(3\om_1)}(wq^3t^{-1}) -
  \de( q^{-6}t^{2}{w\over z}) T_{V(3\om_1)}(zq^3t^{-1}) \biggr) \cr
  & + \cst{q^3,q^3,t^{-1},qt^{-1},q^5t^{-1},q^2t^{-2}}
  {q,q^2t^{-1},q^3t^{-1},q^3t^{-1},q^5t^{-2}} \cr
  & \times
  \biggl( \de( q^{10}t^{-4}{w\over z}) T_{V'(2\om_1)}(wq^5t^{-2}) -
  \de( q^{-10}t^{4}{w\over z}) T_{V'(2\om_1)}(zq^5t^{-2}) \biggr) \cr
 & + \cst{qt,q^2,q^3,t^{-1},q^4t^{-1},qt^{-2}}
  {q,qt^{-1},q^2t^{-1},q^3t^{-1},q^4t^{-2}} \cr
 & \times
 \biggl( \de( q^8t^{-4}{w\over z}) T_2(wq^4t^{-2} ) -
 \de( q^{-8}t^{4} {w\over z}) T_2(z q^4t^{-2}) \biggr) \cr
 & + \cst{q^3,t^{-1},q^4t^{-1},q^5t^{-1},qt^{-2},q^2t^{-2},q^6t^{-2},q^3t^{-3}}
 {qt^{-1},q^2t^{-1},q^3t^{-1},q^3t^{-2},q^4t^{-2},q^5t^{-2},q^6t^{-3}}
 \biggl( \de( q^{12}t^{-6}{w\over z}) - \de( q^{-12}t^{6}{w\over z}) \biggr) 
  \,, \cr}
}
where 
\eqn\eqBGb{ \eqalign{
T_{V(3\om_1)}(z) & \eql \NO{ Y_1(zq^{-2})Y_1(z)Y_1(zq^2) } + \ldots\,,\cr
T_{V(2\om_1)}(z) & \eql \NO{ Y_1(zq^{-1})Y_1(zq) } + \ldots\,,\cr
T_{V'(2\om_1)}(z) & \eql \NO{ Y_1(zq^{-4}t)Y_1(zq^4t^{-1}) } + \ldots\,,\cr}
}
are $\cW_{q,t}[G_2]$ generators corresponding to irreducible 
$U_q(\widehat{G_2})$ representations that decompose under $U_q(G_2)$ as
${\bf 77'}\oplus {\bf 27} \oplus 2 ({\bf 14})$,
${\bf 27}\oplus {\bf 7}$ and ${\bf 27}\oplus {\bf 14} \oplus {\bf 1}$,
respectively.


\footatend\immediate\closeout\rfile
\baselineskip=14pt{\bigskip\noindent {\bf References}}%
\bigskip{\frenchspacing%
\parindent=20pt\escapechar=` \input refs.tmp\vfill\eject}\nonfrenchspacing
\vfill\eject\end